\documentclass[11pt, a4paper]{article}
\usepackage{amsmath,amssymb,latexsym,color}
\usepackage[mathscr]{eucal}
\usepackage[colorlinks,
            linkcolor=blue,
            anchorcolor=green,
            citecolor=magenta
           ]{hyperref}
\newtheorem{thm}{Theorem}[section]

\newtheorem{lemma}[thm]{Lemma}

\newtheorem{rem}{Remark}

\newtheorem{ass}{Assumption}

\newtheorem{example}{Example}[section]
\newtheorem{defin}{Definition}[section]

\newcommand{\proof}{{\it Proof.\quad}}
\newcommand{\qed}{\hfill$\Box$\medskip}

\usepackage{CJK}

\usepackage{authblk}

\long\def\delete#1{}

\usepackage{color}

\definecolor{Blue}{rgb}{0,0,1}
\definecolor{Red}{rgb}{1,0,0}
\definecolor{DarkGreen}{rgb}{0,0.6,0}
\definecolor{DarkYellow}{rgb}{1,1,0.2}
\definecolor{DarkPurple}{rgb}{.6,0,1}

\usepackage{xcolor}
\usepackage[normalem]{ulem}

\usepackage{cleveref}
\crefformat{section}{\S#2#1#3}
\crefformat{subsection}{\S#2#1#3}
\crefformat{subsubsection}{\S#2#1#3}
\crefrangeformat{section}{\S\S#3#1#4 to~#5#2#6}
\crefmultiformat{section}{\S\S#2#1#3}{ and~#2#1#3}{, #2#1#3}{ and~#2#1#3}

\begin{document}
%\begin{CJK*}{GBK}{song}

\renewcommand{\baselinestretch}{1.3}
%%%%%%%%%%%%%%%%%%%%%%%%%%%%%%%%%%%%%%%%%%%%%%%%%%%%%%%%%%%%%%%%%%%%%%%%%%%%%%%%%%%%%%%%
%%%%%%%%%%%%%%%%%%%%%%%%%%%%%%%%%%%%%%%%%%%%%%%%%%%%%%%%%%%%%%%%%%%%%%%%%%%%%%%%%%%%%%%%
\title{\bf Non-trivial $t$-intersecting families for vector spaces}

\author[1,2]{Mengyu Cao\thanks{E-mail: \texttt{caomengyu@mail.bnu.edu.cn}}}
\author[1]{Benjian Lv\thanks{Corresponding author. E-mail: \texttt{bjlv@bnu.edu.cn}}}
\author[1]{Kaishun Wang\thanks{E-mail: \texttt{wangks@bnu.edu.cn}}}
\author[3]{Sanming Zhou\thanks{E-mail: \texttt{sanming@unimelb.edu.au}}}
\affil[1]{\small Laboratory of Mathematics and Complex Systems (Ministry of Education), School of Mathematical Sciences, Beijing Normal University, Beijing 100875, China}
\affil[2]{\small Department of Mathematical Sciences, Tsinghua University, Beijing 100084, China}
\affil[3]{\small School of Mathematics and Statistics, The University of Melbourne, Parkville, VIC 3010, Australia}

\date{}
\maketitle

\begin{abstract}
Let $V$ be an $n$-dimensional vector space over a finite field $\mathbb{F}_q$.
In this paper we describe the structure of maximal non-trivial $t$-intersecting families of $k$-dimensional subspaces of $V$ with large size.
We also determine the non-trivial $t$-intersecting families with maximum size. In the special case when $t=1$ our result gives rise to the well-known Hilton-Milner Theorem for vector spaces.

\medskip
\noindent {\em AMS Classification (2020):} 05D05, 05A30


\end{abstract}

\section{Introduction}

The study of intersecting families has long been an important area of research in combinatorics \cite{Deza-Frankl-1983, GK} ever since the birth of the celebrated Erd\H{o}s-Ko-Rado Theorem \cite{Erdos-Ko-Rado-1961-313}. In this paper we give a description of the structure of maximal non-trivial $t$-intersecting families of $k$-subspaces of an $n$-dimensional vector space over a finite field whose size is a bit smaller than the bound in the Erd\H{o}s-Ko-Rado Theorem for vector spaces. In particular, we extend the Hilton-Milner Theorem for vector spaces \cite{AB} by describing the structure of non-trivial $t$-intersecting families of vector spaces with maximum size.

Let $n$ and $k$ be integers with $1\leq k\leq n.$ Write $[n]=\{1,2,\ldots,n\}$ and denote by ${[n]\choose k}$ the family of all $k$-subsets of $[n].$ For any positive integer $t$, a family $\mathcal{F}\subseteq {[n]\choose k}$ is said to be \emph{$t$-intersecting} if $|A \cap B|\geq t$ for all $A, B\in\mathcal{F}.$ A family is called \emph{intersecting} if it is $1$-intersecting. A $t$-intersecting family is called \emph{trivial} if all its members contain a common specified $t$-subset of $[n]$, and \emph{non-trivial} otherwise.

The Erd\H{o}s-Ko-Rado Theorem gives the maximum size of a $t$-intersecting family and shows further that any $t$-intersecting family with maximum size is a trivial family consisting of all $k$-subsets that contain a fixed $t$-subset of $[n]$ for $n>(t +1)(k-t +1)$ \cite{Erdos-Ko-Rado-1961-313,Frankl-1978,Wilson-1984}. In \cite{Ahlswede-Khachatrian-1997,Frankl--Furedi-1991}, the structure of such extremal families for any positive integers $t,k$ and $n$ was described.
Determining the structure of non-trivial $t$-intersecting families of $k$-subsets of $[n]$ with maximum size was a long-standing problem. The first such result is the Hilton-Milner Theorem \cite{Hilton-Milner-1967} which describes the structure of such families for $t=1$. A complete solution to this problem for any $t$ was obtained by Ahlswede and Khachatrian \cite{Ahlswede-Khachatrian-1996}.
Recently, other maximal non-trivial intersecting families with large size have been studied. For example, Kostochka and Mubayi \cite{Kostochka-Mubayi} described the structure of intersecting families of $k$-subsets of $[n]$ whose size is quite a bit smaller than the bound ${n-1\choose k-1}$ given by the Erd\H{o}s-Ko-Rado Theorem. In \cite{Han-Kohayakawa}, Han and Kohayakawa determined the maximum size of an intersecting family which is not a subfamily of any largest or second largest maximal intersecting family, and characterized all families achieving that extremal value.

The Erd\H{o}s-Ko-Rado Theorem and the Hilton-Milner Theorem for finite sets have natural extensions to vector spaces. Let $n$ and $k$ be integers with $1 \le k \le n$, and $V$ an $n$-dimensional vector space over the finite field $\mathbb{F}_q$, where $q$ is necessarily a prime power. We use ${V\brack k}$ to denote the family of all $k$-dimensional subspaces of $V$.  In the sequel we will abbreviate ``$k$-dimensional subspace'' to ``$k$-subspace''. Recall that  for any positive integers $a$ and $b$ the \emph{Gaussian binomial coefficient} is defined by
$$
{a\brack b} = \prod_{0\leq i<b}\frac{q^{a-i}-1}{q^{b-i}-1}.
$$
In addition, we set ${a\brack 0}=1$ and ${a\brack c} =0$ if $c$ is a negative integer. It is well known that the size of ${V\brack k}$ is equal to ${n\brack k}$.

For any positive integer $t$, a family $\mathcal{F}\subseteq{V\brack k}$ is called $t$-\emph{intersecting} if $\dim(A\cap B)\geq t$ for all $A,B\in\mathcal{F}$. A family is called \emph{intersecting} if it is $1$-intersecting. A $t$-intersecting family $\mathcal{F}\subseteq{V\brack k}$ is called \emph{trivial} if all its members contain a common specified $t$-subspace of $V$ and \emph{non-trivial} otherwise. In general, the triviality of an intersecting family is determined by the following parameter introduced in \cite{AB}:
For any $\mathcal{F}\subseteq{V\brack k}$, the \emph{covering number} $\tau(\mathcal{F})$ of $\mathcal{F}$ is the minimum dimension of a subspace $T$ of $V$ such that $\dim(T\cap F)\geq 1$ for every $F\in \mathcal{F}$. It is clear that an intersecting family $\mathcal{F}$ is trivial if and only if $\tau(\mathcal{F})=1$.

Let $n,k$ and $t$ be positive integers with $n \ge 2k\ge 2t$, and $\mathcal{F}\subseteq{V\brack k}$ a $t$-intersecting family with maximum size. The Erd\H{o}s-Ko-Rado Theorem for vector spaces shows that $\mathcal{F}$ must be a trivial family consisting of all $k$-subspaces of $V$ which contain a fixed $t$-subspace of $V$, or $n = 2k$ and $\mathcal{F}$ consists of all $k$-subspaces of a fixed $(n-t)$-subspace of $V$ \cite{Deza-Frankl-1983,PR,Hsieh-1975-1,Tanaka-2006-903}. Using the covering number, Blokhuis et al. \cite{AB} obtained a vector space version of the Hilton-Milner Theorem, which described the structure of any non-trivial intersecting family with maximum size.

In this paper we study maximal non-trivial $t$-intersecting families of $k$-subspaces of $V$ for any positive integer $t$. By \cite[Remark (ii) in Section 9.3]{ABE} any maximal non-trivial $(k-1)$-intersecting family of $k$-subspaces of $V$ is the collection of all $k$-subspaces contained in a fixed $(k+1)$-subspace of $V$. Henceforth we will only consider the case when $6\leq 2k\leq n$ and $1\leq t\leq k-2$. 

To present our results let us first introduce the following three constructions of $t$-intersecting families of $k$-subspaces of $V$.
\medskip

\noindent\textbf{Family I.}\quad Let $X$ and $M$ be subspaces of $V$ such that $X\subseteq M$, $\dim(X)=t$ and $\dim(M)=k+1.$ Define
\begin{equation}
\label{eq:H1}
\mathcal{H}_1(X,M)=\left\{F\in{V\brack k}\mid X\subseteq F,\ \dim(F\cap M)\geq t+1\right\}\cup{M\brack k}.
\end{equation}
%\begin{equation}
%\label{eq:h1}
%h_1(t,k+1) = |\mathcal{H}_1(X,M)|.
%\end{equation}

\noindent\textbf{Family II.}\quad Let $X,M$ and $C$ be subspaces of $V$ such that $X\subseteq M\subseteq C,$ $\dim(X)=t$, $\dim(M)=k$ and $\dim(C)=c$, where $c\in\{k+1,k+2,\ldots, 2k-t,n\}$. Define
\begin{align}
\mathcal{H}_2(X,M,C) &= \mathcal{A}(X,M) \cup \mathcal{B}(X, M,C) \cup \mathcal{C}(X, M, C)\label{eq:H2},
\end{align}
where
\begin{eqnarray*}
\mathcal{A}(X,M) &=&\left\{F\in{V\brack k}\mid X\subseteq F,\ \dim(F\cap M)\geq t+1\right\},\\% \label{eq:A}\\
\mathcal{B}(X, M,C) &=&\left\{F\in{V\brack k}\mid F\cap M=X,\ \dim(F\cap C)=c-k+t\right\},\\% \label{eq:B}\\
\mathcal{C}(X, M, C) &=&\left\{F\in{C\brack k}\mid \dim(F\cap X)=t-1,\ \dim(F\cap M)=k-1\right\}. %\label{eq:C}
\end{eqnarray*}

\noindent\textbf{Family III.}\quad Let $Z$ be a $(t+2)$-subspace of $V$. Define
\begin{equation}
\label{eq:H4}
\mathcal{H}_3(Z)=\left\{F\in{V\brack k}\mid \dim(F\cap Z)\geq t+1\right\}.
\end{equation}

It is straightforward to verify that $\mathcal{H}_1(X,M), \mathcal{H}_2(X,M,C), \mathcal{H}_2(X,M,V)$ and $\mathcal{H}_3(Z)$ are all non-trivial $t$-intersecting families of $k$-subspaces of $V$.

\begin{rem}
\label{rem1}
{\em
In Family II, if $C$ satisfies $\dim(C) = k+1$, then $\mathcal{H}_2(X,M,C) = \mathcal{H}_1(X,C)$; if $t$ and $k$ satisfy $t=k-2$, then $\mathcal{H}_2(X,M,V)=\mathcal{H}_3(M).$
}
\end{rem}
\iffalse
Define
\begin{equation}
\label{eq:f}
f(n,k,t)={k-t\brack 1}{n-t-1\brack k-t-1}-q{k-t\brack 2}{n-t-2\brack k-t-2}.
\end{equation}
\fi

Our first main result describes the structure of all maximal non-trivial $t$-intersecting families of $k$-subspaces of $V$ with large size.

\begin{thm}
\label{main-1}
Let $n, k$ and $t$ be positive integers with $t\leq k-2$ and $2k+t+\min\{4, 2t\}\leq n$. If $\mathcal{F}\subseteq {V\brack k}$ is a maximal non-trivial $t$-intersecting family and
$$|\mathcal{F}| \geq {k-t\brack 1}{n-t-1\brack k-t-1}-q{k-t\brack 2}{n-t-2\brack k-t-2},$$ then one of the following holds:
\begin{itemize}
\item[{\rm(i)}] $\mathcal{F}=\mathcal{H}_2(X,M,C)$ for some $t$-subspace $X$, $k$-subspace $M$ and $c$-subspace $C$ of $V$ with $X\subseteq M\subseteq C$ and $c\in\{k+1,k+2,\ldots, 2k-t,n\}$;
\item[{\rm(ii)}] $\mathcal{F}=\mathcal{H}_3(Z)$ for some $(t+2)$-subspace $Z$ of $V$, and $\frac{k}{2}-1\leq t\leq k-2$.
\end{itemize}
\end{thm}

By comparing the size of the families given in Theorem \ref{main-1}, we can describe the structure of the non-trivial $t$-intersecting families with maximum size. Our second main result is as follows.

\begin{thm}
\label{main-2}
Let $n, k$ and $t$ be positive integers with $t\leq k-2$ and $2k+t+\min\{4, 2t\}\leq n$. Then, for any non-trivial $t$-intersecting family $\mathcal{F}\subseteq {V\brack k}$, the following hold:
\begin{itemize}
\item[{\rm(i)}]  if $1\leq t\leq \frac{k}{2}-1$, then
$$
|\mathcal{F}|\leq {n-t\brack k-t}-q^{(k+1-t)(k-t)}{n-k-1\brack k-t}+q^{k+1-t}{t\brack 1},
$$
and equality holds if and only if $\mathcal{F}=\mathcal{H}_1(X,M)$ for some $t$-subspace $X$ and $(k+1)$-subspace $M$ of $V$ with $X\subset M$;
\item[{\rm(ii)}] if $\frac{k}{2}-1<t\leq k-2,$ then
$$
|\mathcal{F}|\leq {t+2\brack 1}{n-t-1\brack k-t-1}-q{t+1\brack 1}{n-t-2\brack k-t-2},
$$
and equality holds if and only if $\mathcal{F}=\mathcal{H}_3(Z)$ for some $(t+2)$-subspace $Z$ of $V$, or $(t, k) = (1, 3)$ and $\mathcal{F}=\mathcal{H}_1(X,M)$ for some $1$-subspace $X$ and $3$-subspace $M$ of $V$ with $X\subset M$.
\end{itemize}
\end{thm}

In the special case when $t=1$, Theorem~\ref{main-2}  gives rise to the Hilton-Milner Theorem for vector spaces with $n\geq 2k+3$ (\cite{AB}).

The rest of this paper is organized as follows. In the next section we will prove a number of inequalities for the sizes of the intersecting families in Families~I, II and III. In \cref{sec:ub} we will prove some upper bounds for the sizes of non-trivial $t$-intersecting families of subspaces of $V$ using a key notion---$t$-covering number, which is a generalization of the covering number. After these preparations we will prove Theorems \ref{main-1} and \ref{main-2} in \cref{sec:pf}.

\section{Inequalities for the sizes of the constructed families}

\subsection{Equalities and formulas involving the Gaussian binomial coefficients}

This subsection is a prepration for \cref{subsec:ineq} and \cref{sec:ub}. The following lemma can be easily proved.

\begin{lemma}\label{lem1-1-1}
Let $m$ and $i$ be positive integers with $i\leq m.$ Then the following hold:
\begin{itemize}
\item[{\rm (i)}] ${m\brack i}={m-1\brack i-1}+q^i{m-1\brack i}$ and ${m\brack i}=\frac{q^m-1}{q^i-1}\cdot{m-1\brack i-1}$;
\item[{\rm (ii)}] $q^{m-i}<\frac{q^m-1}{q^i-1}<q^{m-i+1}$ and $q^{i-m-1}<\frac{q^i-1}{q^m-1}<q^{i-m}$ if $i <m$;
\item[{\rm(iii)}] $q^{i(m-i)}\leq{m\brack i}< q^{i(m-i+1)}$, and $q^{i(m-i)}<{m\brack i}$ if $i <m$;
\item[{\rm (iv)}] $\frac{q^m-1}{q^i-1}<2q^{m-i}$.
\end{itemize}
\end{lemma}

\iffalse
\begin{lemma}{\rm{(\cite[Proposition 2.3]{KWW})}}\label{2.0}
Let $W$ be an $n$-dimensional vector space over $\mathbb{F}_q$. For $1\leq m\leq n$ and $0\leq i\leq \min
\{m,n-m\}$, let $P'$ and $Q'$ be two fixed $m$-dimensional subspaces of $W$ with $\dim(P'\cap Q')=m-i$. Then the number of $m$-dimensional subspaces $S'$ of $W$ satisfying $\dim(P'\cap S')=m-s$ and $\dim(S'\cap Q')=m-u$ is
\begin{eqnarray*}
p^{i}_{s,u}(m,n)=&\sum_{\rho+\alpha=u,\beta+\gamma=m-u,\rho+\gamma\leq s}q^{\omega}\prod_{k=i+\rho-s+1}^{i-\gamma}(q^k-1)\\
&\times\textstyle{{\alpha\brack s-\rho-\gamma}{i\brack \alpha}{i\brack \gamma}{m-i\brack \beta}{n-m-i\brack \rho}},
\end{eqnarray*}
where $\omega=\frac{1}{2}(s-\gamma-\rho)(s-\gamma-\rho-1)+(m-\beta)(m-\beta-i)+\rho(2i-\alpha-\gamma)$.
\end{lemma}
\fi

Set
\begin{eqnarray*}
g_1(t,n)&=&{t+2\brack 1}{n-t-1\brack t+1}-q{t+1\brack 1}{n-t-2\brack t}\\
g_2(t,n)&=&{n-t\brack t+2}-q^{(t+2)^2}{n-2t-2\brack t+2}.
\end{eqnarray*}

\begin{lemma}\label{lem9}
We have
$$
g_1(t,n)-g_2(t,n)=\sum_{j=1}^{t}q^{j(t+2)+1}{t+1-j\brack 1}{n-t-2-j\brack t}.
$$
\end{lemma}
\proof
By Lemma~\ref{lem1-1-1}(i), we have
\begin{align*}
g_1(t,n) &= {n-t-1\brack t+1}+q{t+1\brack 1}{n-t-1\brack t+1}-q{t+1\brack 1}{n-t-2\brack t}\\
&= {n-t-1\brack t+1}+q^{t+2}{t+1\brack 1}{n-t-2\brack t+1}.
\end{align*}
Using Lemma~\ref{lem1-1-1}(i) repeatedly, we can show that
$$
g_2(t,n) = \sum_{i=1}^{t+2}q^{(t+2)(i-1)}{n-t-i\brack t+1}.
$$
Set
\begin{eqnarray*}
f(a) &=&\ \sum_{j=1}^{a}q^{j(t+2)+1}{t+1-j\brack 1}{n-t-2-j\brack t}+q^{(a+1)(t+2)}{t+1-a\brack 1}{n-t-2-a\brack t+1}\\
& & -\sum_{i=a+2}^{t+2}q^{(t+2)(i-1)}{n-t-i\brack t+1}
\end{eqnarray*}
for $a \in \{0,1,\ldots,t\}$.
Then
$$
g_1(t,n) - g_2(t,n) = q^{t+2}{t+1\brack 1}{n-t-2\brack t+1}-\sum_{i=2}^{t+2}q^{(t+2)(i-1)}{n-t-i\brack t+1} = f(0).
$$
On the other hand, by Lemma~\ref{lem1-1-1}(i), we have
\begin{eqnarray*}
f(a+1)-f(a) & = & q^{(a+1)(t+2)+1}{t-a\brack 1}{n-t-3-a\brack t}+q^{(a+2)(t+2)}{t-a\brack 1}{n-t-3-a\brack t+1}\\
& & -q^{(a+1)(t+2)}{t+1-a\brack 1}{n-t-2-a\brack t+1}+q^{(t+2)(a+1)}{n-t-a-2\brack t+1} \\
& = & 0.
\end{eqnarray*}
Since this holds for each $a$, we obtain $f(0)=f(1)=\cdots=f(t)$. Therefore,
$$
g_1(t,n)-g_2(t,n) = f(0) = f(t) = \sum_{j=1}^{t}q^{j(t+2)+1}{t+1-j\brack 1}{n-t-2-j\brack t}
$$
as required. \qed

Let $W$ be an $(e+l)$-dimensional vector space over $\mathbb{F}_q$, where $l,e \ge 1$, and let $L$ be a fixed $l$-subspace of $W$. We say that an $m$-subspace $U$ is of \emph{type} $(m,h)$ if $\dim(U\cap L)=h$. Define $\mathcal{M}(m,h;e+l,e)$ to be the set of all subspaces of $W$ with type $(m,h)$. 
\begin{lemma}{\rm(\cite[Lemma~2.1]{KWW})}\label{lem4}
$\mathcal{M}(m,k;e+l,e)$ is non-empty if and if $0\leq h\leq l$ and $0\leq m-h\leq e.$ Moreover, if $\mathcal{M}(m,h;e+l,e)$ is non-empty, then
$$
|\mathcal{M}(m,h;e+l,e)|=q^{(m-h)(l-h)}{e\brack m-h}{l\brack h}.	
$$
\end{lemma}

Define
$$
N'(m_{1},h_{1};m,h;e+l,e)
$$
to be the number of subspaces of $W$ with type $(m,h)$ containing a given subspace with type $(m_{1},h_{1})$. Observe that $|\mathcal{M}(m,h;e+l,e)| = N'(0,0;m,h;e+l,e)$.

\begin{lemma}{\rm{(\cite{KW})}}\label{lem5}
$N^{'}(m_{1},h_{1};m,h;e+l,e)\not= 0$ if and only if $0\leq h_{1}\leq h\leq l$ and $0\leq m_{1}-h_{1}\leq m-h\leq e$. Moreover, if $N^{'}(m_{1},h_{1};m,h;e+l,e)\neq 0,$ then
 $$ N^{'}(m_{1},h_{1};m,h;e+l,e)=q^{(l-h)(m-h-m_{1}+h_{1})}{{e-(m_{1}-h_{1})}\brack{(m-h)-(m_{1}-h_{1})}}{{l-h_{1}}\brack{h-h_{1}}}.
 $$
\end{lemma}

Let 
$$
h_1(t,k+1)=|\mathcal{H}_1(X,M)|,
$$ 
$$
h_2(t,k,c)=|\mathcal{H}_2(X,M,C)|,\ \mbox{for}\ c\in\{k+1,k+2,\ldots,2k-t,n\},
$$
and 
$$
h_3(t+2)=|\mathcal{H}_3(Z)|.
$$ 
The following lemma gives the sizes of  Families I, II and III. 

\begin{lemma}\label{lem6}
Suppose $c\in\{k+1,k+2,\ldots,2k-t,n\}$. Then the following hold:
\begin{align}
h_1(t,k+1)=&\ {n-t\brack k-t}-q^{(k+1-t)(k-t)}{n-k-1\brack k-t}+q^{k+1-t}{t\brack 1}; \label{hmtk+1}\\
h_2(t,k,c)=&\ {n-t\brack k-t}-q^{(k-t)^2}{n-k\brack k-t}+q^{(k-t)^2}{n-c\brack 2k-c-t}+q^{k-t+1}{c-k\brack 1}{t\brack 1}\label{hmtkc};\\
h_3(t+2)=&\ {t+2\brack 1}{n-t-1\brack k-t-1}-q{t+1\brack 1}{n-t-2\brack k-t-2}.\label{hmt+2}
\end{align}
\end{lemma}
\proof Suppose that $X,\ M$ and $C$ are subspaces of $V$ with $X\subseteq M\subseteq C$ such that $\dim(X)=t$, $\dim(M)=k$ and $\dim(C)=c$.  %Let $\mathcal{A} = \mathcal{A}(X, M),\ \mathcal{B} = \mathcal{B}(X, M, C)$ and $\mathcal{C} = \mathcal{C}(X,M,C)$ be the families defined in Construction~\ref{con2}.
Then
\begin{equation*}
\mathcal{A}(X, M)= \left\{F\in{V\brack k}\mid X\subseteq F\right\}\setminus\left\{F\in{V\brack k}\mid F\cap M=X\right\},
%\mathcal{B}&=&\bigcup_{D\subseteq C, D\cap M=X,\atop \dim(D)=c-k+t}\left\{F\in{V\brack k}\mid F\cap C=D\right\},\\
%\mathcal{C}&=&\left\{F\in{C\brack k}\mid \dim(F\cap M)=k-1\right\}\setminus\left\{F\in{C\brack k}\mid \dim(F\cap M)=k-1,\ X\subseteq F\right\}.
\end{equation*}
which implies that $|\mathcal{A}(X, M)|={n-t\brack k-t}- q^{(k-t)^2}{n-k\brack k-t}$ by Lemmas~\ref{lem4} and \ref{lem5}. From Corollary 2.3 and Lemma 2.1 in \cite{Guo}, we have
\begin{equation*}
%|\mathcal{A}|&=&\ {n-t\brack k-t}-N^\prime(t,t;k,t;n,n-k)\\
%&=& {n-t\brack k-t}- q^{(k-t)^2}{n-k\brack k-t},\\
|\mathcal{B}(X, M, C)|= q^{(k-t)^2}{n-c\brack 2k-c-t}\quad \mbox{and}\quad |\mathcal{C}(X,M,C)|= q^{k-t+1}{c-k\brack 1}{t\brack 1}.
\end{equation*}
Since $h_{2}(t, k, c) = |\mathcal{A}(X, M)|+|\mathcal{B}(X, M, C)|+|\mathcal{C}(X,M,C)|$, we obtain (\ref{hmtkc}) immediately.  By Remark~\ref{rem1} and Lemma~\ref{lem1-1-1}(i), we obtain (\ref{hmtk+1}).

Consider the family $\mathcal{H}_3(Z)$, where $Z$ is a $(t+2)$-subspace of $V$. By Lemma~\ref{lem4}, the number of $k$-subspaces $F$ of $V$ satisfying $\dim(F\cap Z)=t+1$ is $q^{k-t-1}{n-t-2\brack k-t-1}{t+2\brack 1}$, and the number of $k$-subspaces $F$ of $V$ satisfying $\dim(F\cap Z)=t+2$ is ${n-t-2\brack k-t-2}.$ Combining these with Lemma~\ref{lem1-1-1}(i), we obtain (\ref{hmt+2}) immediately.
\qed

\subsection{Inequalities for $h_1(t,k+1), h_2(t,k,c)$ and $h_3(t+2)$}
\label{subsec:ineq}

\begin{lemma}\label{lem7}
Let $n, k$ and $t$ be positive integers with $6\leq 2k\leq n$ and $1\leq t\leq k-2.$
\begin{itemize}
\item[{\rm(i)}] We have
$$
h_1(t,k+1) = h_2(t,k,k+1),
$$
and
$$
h_2(t,k,c)>h_2(t,k,c+1)
$$
for $c\in\{k+1,k+2,\ldots,2k-t-1\}$.
\item[{\rm(ii)}] Assume that $1\leq t\leq k-3$. If $2k\leq n\leq (k-t)^2-1$, or $n=(k-t)^2$ and $q\geq 3$, or $(n,q,t)=((k-t)^2,2,1),$ then
$$
h_2(t,k,2k-t)>h_2(t,k,n).
$$
If $n\geq (k-t)^2+1$, or $(n,q)=((k-t)^2,2)$ and $t\geq 2,$ then
$$
h_2(t,k,k+1)>h_2(t,k,n)>h_2(t,k,2k-t).
$$
\item[{\rm(iii)}] Assume that $t=k-2$. If $t=1$, then
$$
h_2(t,k,n)=h_2(t,k,k+1);
$$
and if $t \ge 2$, then
$$
h_2(t, t+2, n)>h_2(t, t+2, t+3).
$$
\end{itemize}
\end{lemma}
\proof (i)\quad As seen in Remark~\ref{rem1} we have $h_1(t,k+1)=h_2(t,k,k+1)$ for $1\leq t\leq k-2$. For $c\in\{k+1,k+2,\ldots,2k-t-1\}$, by Lemma~\ref{lem1-1-1}(i), we have
$$
h_2(t,k,c)-h_2(t,k,c+1)=q^{(k-t)^2+2k-c-t}{n-c-1\brack 2k-c-t}-q^{c-t+1}{t\brack 1}.
$$
Since ${n-c-1\brack 2k-c-t}\geq q^{(2k-c-t)(n-2k+t-1)}$ and ${t\brack 1}<q^t$ by Lemma~\ref{lem1-1-1}(iii), and since
\begin{eqnarray*}
& &(k-t)^2+2k-c-t+(2k-c-t)(n-2k+t-1)-c-1\\
& = &(k-t-1)^2+(2k-c-t-1)(n-2k+t+1)+(n-2k)-1\\
& \geq & 0,
\end{eqnarray*}
we obtain $h_2(t,k,c)-h_2(t,k,c+1)>0$ for $c\in\{k+1,k+2,\ldots,2k-t-1\}$.

(ii) \quad Note that
\begin{eqnarray}\label{equ13}
h_2(t,k,n)-h_2(t,k,c)=q^{c-t+1}{t\brack 1}{n-c\brack 1}-q^{(k-t)^2}{n-c\brack 2k-c-t}
\end{eqnarray}
for any $c\in\{k+1,k+2,\ldots,2k-t\}$. When $c=k+1$, by (\ref{equ13}) and Lemma~\ref{lem1-1-1}(iii), and noting that $t\leq k-3$ and $2k\leq n$, we obtain
\begin{align*}
h_2(t,k,n)-h_2(t,k,k+1)=&\ q^{k-t+2}{t\brack 1}{n-k-1\brack 1}-q^{(k-t)^2}{n-k-1\brack k-t-1}\\
<&\ q^{n+1}-q^{(n-k+1)(k-t-1)+1}\\
<&\ 0.
\end{align*}

When $c=2k-t$, by (\ref{equ13}) again, we have
\begin{eqnarray*}
h_2(t,k,n)-h_2(t,k,2k-t)=q^{2k-2t+1}{t\brack 1}{n-2k+t\brack 1}-q^{(k-t)^2}.
\end{eqnarray*}
By Lemma~\ref{lem1-1-1}(iii), if $2k\leq n\leq (k-t)^2-1$, then
\begin{eqnarray*}
h_2(t,k,n)-h_2(t,k,2k-t)< q^{n+1}-q^{(k-t)^2}\leq 0;
\end{eqnarray*}
and if $n\geq (k-t)^2+1$, then
\begin{eqnarray*}
h_2(t,k,n)-h_2(t,k,2k-t)> q^{n-1}-q^{(k-t)^2}\geq 0.
\end{eqnarray*}

Now assume that $n=(k-t)^2$. Then
\begin{align*}
h_2(t,k,n)-h_2(t,k,2k-t)=&\ \frac{1}{(q-1)^2}\left((q^t-1)(q^{n-t+1}-q^{2k-2t+1})-q^{(k-t)^2}(q-1)^2\right)\\
=&\ \frac{1}{(q-1)^2}\left(q^{(k-t)^2}(-q^2+3q-1-q^{-t+1})-q^{2k-2t+1}(q^t-1)\right).
\end{align*}
If $q\geq 3$, then $-q^2+3q\leq0,$ and $h_2(t,k,n)-h_2(t,k,2k-t)<0.$
If $q=2$ and $t=1$, then
$h_2(t,k,n)-h_2(t,k,2k-t)=-2^{2k-1}<0.$
If $q=2$ and $t\geq 2$, then
\begin{eqnarray*}
h_2(t,k,n)-h_2(t,k,2k-t)&=&2^{(k-t)^2}(1-2^{-t+1})-2^{2k-2t+1}(2^t-1)\\
&\geq& 2^{(k-t)^2-1}-2^{2k-t+1}+2^{2k-2t+1}\\
& > & 0
\end{eqnarray*}
as $(k-t)^2 = n \geq 2k$.

(iii)\quad If $t=k-2$, then by (\ref{equ13}) we have
$$
h_2(t,k,n)-h_2(t,k,k+1)=q^4{k-2\brack 1}{n-k-1\brack 1}-q^4{n-k-1\brack 1}.
$$
It is clear that $h_2(t,k,n)-h_2(t,k,k+1)=0$ if $t=1$ and $h_2(t,k,n)-h_2(t,k,k+1)>0$ if $t\geq2$.
\qed

Define
\begin{equation}
\label{eq:f}
f(n,k,t)={k-t\brack 1}{n-t-1\brack k-t-1}-q{k-t\brack 2}{n-t-2\brack k-t-2}.
\end{equation}

\begin{lemma}\label{lem8}
Let $n, k$ and $t$ be positive integers with $6\leq 2k\leq n$ and $1\leq t\leq k-2$.% Let $f(n,k,t)$ be the function defined in \eqref{eq:f}.
\begin{itemize}
\item[{\rm(i)}] $\min\{h_2(t,k,2k-t),\ h_2(t,k,n)\}\geq f(n,k,t).$
\item[{\rm(ii)}]  If $1\leq t\leq k-3$, then
\begin{eqnarray}\label{ubound-1}
h_1(t,k+1)\leq {k-t+1\brack 1}{n-t-1\brack k-t-1}.
\end{eqnarray}
\item[{\rm(iii)}]
If $1\leq t\leq k-4,$ then
\begin{eqnarray}\label{ubound-2}
h_1(t,k+1)\leq {k-t+1\brack 1}{n-t-1\brack k-t-1}-q^{(k-t-1)(k-t-2)+1}{n-k-1\brack k-t-2}{k+1-t\brack 2}.
\end{eqnarray}
\end{itemize}
\end{lemma}
\proof Let $X$ and $M$ be subspaces of $V$ with $\dim(X)=t$ and $X\subseteq M.$ For each $i\in\{t,t+1,\ldots,k\}$, set
$$
\mathcal{A}_i(X,M)=\left\{F\in{V\brack k}\mid X\subseteq F,\ \dim(F\cap M)=i\right\}
$$
and
$$
\mathcal{L}_i(X,M)=\left\{(I,F)\in{V\brack i}\times{V\brack k}\mid X\subseteq I\subseteq M,\ I\subseteq F \right\}.
$$
Using Lemma~\ref{lem1-1-1} and double counting $|\mathcal{L}_i(X,M)|$, we obtain
\begin{eqnarray}
|\mathcal{L}_i(X,M)|=\sum_{j=i}^k\left|\mathcal{A}_j(X,M)\right|\cdot{j-t\brack i-t}={\dim(M)-t\brack i-t}{n-i\brack k-i}.\label{equ15}
\end{eqnarray}
In particular, we have
\begin{eqnarray}
|\mathcal{L}_{t+1}(X,M)|=\sum_{j=t+1}^k\left|\mathcal{A}_j(X,M)\right|+\sum_{j=t+2}^k\left|\mathcal{A}_j(X,M)\right|\cdot\left({j-t\brack 1}-1\right)\label{equ14}.
\end{eqnarray}

(i)\quad Let $M$ be a $k$-subspace of $V$ and $\mathcal{A}(X, M)$ the family constructed in Family~II. Observe that $\mathcal{A}(X, M) = \cup_{j=t+1}^k\mathcal{A}_j(X,M).$ By (\ref{equ15}) and (\ref{equ14}), we obtain
\begin{eqnarray*}
|\mathcal{L}_{t+1}(X,M)|&=& {k-t\brack 1}{n-t-1\brack k-t-1} \\
& \leq & |\mathcal{A}(X, M)|+\sum_{j=t+2}^k\left|\mathcal{A}_j(X,M)\right| q{j-t\brack 2}\\
& = & |\mathcal{A}(X, M)|+q|\mathcal{L}_{t+2}(X,M)| \\
& = & |\mathcal{A}(X, M)|+q{k-t\brack 2}{n-t-2\brack k-t-2}.
\end{eqnarray*}
That is,
\begin{eqnarray*}
|\mathcal{A}(X, M)|\geq {k-t\brack 1}{n-t-1\brack k-t-1}-q{k-t\brack 2}{n-t-2\brack k-t-2}.
\end{eqnarray*}
We then obtain (i) by the definitions of $\mathcal{H}_2(X,M,C)$ and $\mathcal{H}_2(X,M,V)$ and the proof of Lemma~\ref{lem6}.

(ii)\quad Let $M$ be a $(k+1)$-subspace of $V$ and $\mathcal{A}'(X,M) = \cup_{j=t+1}^k\mathcal{A}_j(X,M).$ By (\ref{equ15}) and (\ref{equ14}) again, we have
\begin{eqnarray*}
|\mathcal{L}_{t+1}(X,M)| & = & {k-t+1\brack 1}{n-t-1\brack k-t-1}\\
& = & |\mathcal{A}'(X,M)|+\sum_{j=t+2}^k\left|\mathcal{A}_j(X,M)\right| \left({j-t\brack 1}-1\right).
\end{eqnarray*}
Since $1\leq t\leq k-3$ and $\left|\mathcal{A}_{k-1}(X,M)\right|=N^\prime(t,t; k,k-1;n,n-k-1),$ by Lemma~\ref{lem1-1-1}(iii) and the assumption $n\geq 2k$, we have
\begin{eqnarray*}
\left|\mathcal{A}_{k-1}(X,M)\right| \left({k-1-t\brack 1}-1\right)
&=&q^2{n-k-1\brack 1}{k+1-t\brack k-1-t} \cdot q{k-t-2\brack 1}\\
&>&q^{n+2k-3t-4}>q^{k+1}\geq q^{k+1-t}{t\brack 1}.
\end{eqnarray*}
Observe that, for any $F\in{M\brack k}$, $F\in\mathcal{A}'(X,M)$ if and only if $X\subseteq F$. Thus, by Lemma~\ref{lem4}, we have
$$
\left|{M\brack k}\setminus \mathcal{A}'(X,M)\right|={k+1\brack k}-{k+1-t\brack k-t}=q^{k+1-t}{t\brack 1}.
$$
So by the construction of $\mathcal{H}_1(X,M)$ we then obtain
\begin{eqnarray*}
h_1(t,k+1)\} & = & |\mathcal{A}'(X,M)|+q^{k+1-t}{t\brack 1}\\
& \leq & |\mathcal{L}_{t+1}(X,M)| \\
& = & {k-t+1\brack 1}{n-t-1\brack k-t-1}
\end{eqnarray*}
as required.

(iii)\quad Continuing our discussion in (ii), if $t\leq k-4$, then $t+2\neq k-1$ and hence
\begin{eqnarray*}
h_1(t,k+1) & = & |\mathcal{A}'(X,M)|+q^{k+1-t}{t\brack 1}\\
& \leq & |\mathcal{L}_{t+1}(X,M)|-|\mathcal{A}_{t+2}(X,M)| \left({2\brack 1}-1\right).
\end{eqnarray*}
This together with $|\mathcal{A}_{t+2}(X,M)|=N^\prime(t,t; k,t+2;n,n-k-1)$ and Lemma~\ref{lem5} yields (\ref{ubound-2}).
\qed

\begin{lemma}\label{lem10}
Let $n, k$ and $t$ be positive integers with $6\leq 2k\leq n$ and $1\leq t\leq k-2$. Let $f(n,k,t)$ be the function defined in \eqref{eq:f}.
\begin{itemize}
\item[{\rm(i)}] If $1\leq t< \frac{k}{2}-1,$ then $h_3(t+2)<f(n,k,t).$
\item[{\rm(ii)}] If $\frac{k}{2}-1\leq t\leq k-2,$ then $h_3(t+2)>f(n,k,t).$
\end{itemize}
\end{lemma}
\proof Let
$$
f_1(n,k,t) = \frac{f(n,k,t)-h_3(t+2)}{{n-t-2\brack k-t-2}}.
$$
By (\ref{hmt+2}), we have
\begin{eqnarray*}
f_1(n,k,t)=\frac{(q^{k-t}-q^{t+2})(q^{n-t-1}-1)}{(q-1)(q^{k-t-1}-1)}+q{t+1\brack 1}-q{k-t\brack 2}.
\end{eqnarray*}

(i)\quad  Suppose that $1\leq t< \frac{k}{2}-1$. Since $k>2t+2$ and $n\geq 2k$, by Lemma~\ref{lem1-1-1}(ii)(iii), we have $f_1(n,k,t)>q^{n-t-1}+q^{t+1}-q^{2k-2t-1}>0$,
which implies $f(n,k,t) > h_3(t+2)$ as required.

(ii)\quad  If $t=\frac{k}{2}-1,$ then
\begin{eqnarray*}
f_1(n,k,t)=q{t+1\brack 1}-q{k-t\brack 2}<0.
\end{eqnarray*}
If $\frac{k}{2}-1< t\leq k-2,$ then by Lemma~\ref{lem1-1-1}(ii)(iii) and the assumption $n\geq 2k$, we have $f_1(n,k,t)<-q^{n-k+t+1}+q^{t+2}-q^{2k-2t-3}<0$. In either case we obtain $f(n,k,t) < h_3(t+2)$ as required. \qed

Combining Lemmas \ref{lem7}(i), \ref{lem8}(i) and \ref{lem10}(i), we obtain that if $1\leq t<\frac{k}{2}-1$ then $\min\{h_2(t,k,c), h_2(t,k,n)\}>h_3(t+2)$ for any $c\in\{k+1, k+2,\ldots,2k-t\}.$ The next lemma gives several inequalities involving $h_2(t,k,c)$, $h_2(t,k,n)$ and $h_3(t+2)$ in the case when $\frac{k}{2}-1\leq t\leq k-2.$

\begin{lemma}
\label{lem8-1}
Let $n, k$ and $t$ be positive integers with $6\leq 2k\leq n$ and $\frac{k}{2}-1 \leq t\leq k-2.$
\begin{itemize}
\item[{\rm(i)}] Assume that $t=\frac{k}{2}-1.$ If $t=1$ and $8\leq n\leq 9,$ then $h_2(t,k,2k-t)> h_3(t+2)$; if $t=1$ and $n\geq 10,$ or $t\geq2,$ then
$h_2(t,k,k+1)> h_3(t+2)>h_2(t,k,2k-t)$; if $t=1$, then $h_2(t,k,n)=h_3(t+2)$; if $t\geq 2$, then $h_2(t,k,k+1)> h_3(t+2)>h_2(t,k,n)$.
\item[{\rm(ii)}] Assume that $\frac{k}{2}-\frac{1}{2}\leq t\leq k-3$. Then $h_3(t+2)>h_2(t,k,k+1)$.
\item[{\rm(iii)}] Assume that $t=k-2.$ If $t=1$, then $h_2(t,k,n)=h_3(t+2)=h_2(t,k,k+1)$; if $t\geq 2,$ then $h_2(t,k,n)=h_3(t+2)>h_2(t,k,k+1)$.
\end{itemize}
\end{lemma}
\proof
(i) \quad Assume that $t=\frac{k}{2}-1.$ By (\ref{hmtkc}), (\ref{hmt+2}) and Lemma~\ref{lem9},  we have
\begin{eqnarray*}
& & h_3(t+2)-h_2(t,k,2k-t)\\
& = & \sum_{j=1}^{t}q^{j(t+2)+1}{t+1-j\brack 1}{n-t-2-j\brack t}-q^{(t+2)^2}-q^{t+3}{t+2\brack 1}{t\brack 1}.
\end{eqnarray*}
It is routine to verify that $h_2(t,k,2k-t)> h_3(t+2)$ when $t=1$ and $8\leq n\leq 9$ and $h_2(t,k,2k-t)< h_3(t+2)$ when $t=1$ and $n=10.$  If $t=1$ and $n\geq 11,$ or $t\geq2$, then
\begin{eqnarray*}
h_3(t+2)-h_2(t,k,2k-t) & > & q^{t(t+2)+1}{n-2t-2\brack t}-q^{(t+2)^2}-q^{t+3}{t+2\brack 1}{t\brack 1}\\
& > & q^{t(n-2t)+1}-q^{(t+2)^2}-q^{3t+5}\\
& > & 0.
\end{eqnarray*}

By (\ref{hmtkc}), (\ref{hmt+2}) and Lemma~\ref{lem9} again,  we have
\begin{eqnarray*}
& & h_3(t+2)-h_2(t,k,n)\\
& = & \sum_{j=1}^{t}q^{j(t+2)+1}{t+1-j\brack 1}{n-t-2-j\brack t}-q^{t+3}{n-2t-2\brack 1}{t\brack 1}.
\end{eqnarray*}
It is straightforward to verify that $h_3(t,k,n)=h_4(t+2)$ when $t=1$. If $t\geq 2$, then
\begin{eqnarray*}
h_3(t+2)-h_2(t,k,n) & > & q^{t(t+2)+1}{n-2t-2\brack t}-q^{t+3}{n-2t-2\brack 1}{t\brack 1}\\
& > & q^{t(n-2t)+1}-q^{n+1}\\
& > & 0.
\end{eqnarray*}

By (\ref{hmtkc}), (\ref{hmt+2}) and Lemma~\ref{lem9} again,  we have
$$
h_2(t,k,k+1)-h_4(t+2)=q^{(t+2)^2}{n-2t-3\brack t+1}-\sum_{j=1}^{t}q^{j(t+2)+1}{t+1-j\brack 1}{n-t-2-j\brack t}+q^{t+3}{t\brack 1}.
$$
Since by Lemma~\ref{lem1-1-1}(ii),
$$
q^{j(t+2)+1}{t+1-j\brack 1}{n-t-2-j\brack t}<q^{nt-2t^2+j+2}
$$
for $j\in\{1,2,\ldots,t\}$, we have
$$
\sum_{j=1}^{t}q^{j(t+2)+1}{t+1-j\brack 1}{n-t-2-j\brack t} < q^{nt-2t^2+2} \sum_{j=1}^t q^j<2 q^{nt-2t^2+2+t}.
$$
Since $q\geq 2$, by Lemma~\ref{lem1-1-1}(iii), we have
$$
q^{(t+2)^2}{n-2t-3\brack t+1}> 2\cdot q^{t^2+4t+3}\cdot q^{(t+1)(n-3t-4)}=2 q^{nt-2t^2+n-3t-1}.
$$
Thus $h_2(t,k,k+1)>h_3(t+2)$ as $n\geq 2k=4t+4.$

(ii)\quad Assume that $\frac{k}{2}-\frac{1}{2}< t\leq k-3.$ By (\ref{hmt+2}) and (\ref{ubound-1}), we have
$$
\frac{h_3(t+2)}{{n-t-2\brack k-t-2}}=\frac{q^{n-t-1}-1}{q^{k-t-1}-1} {t+2\brack 1}-q{t+1\brack 1}
$$
and
$$
\frac{h_1(t,k+1)}{{n-t-2\brack k-t-2}}\leq\frac{q^{n-t-1}-1}{q^{k-t-1}-1}{k-t+1\brack 1}.
$$
Since by Lemma~\ref{lem1-1-1} (ii),
\begin{eqnarray*}
& & \frac{q^{n-t-1}-1}{q^{k-t-1}-1}\cdot{t+2\brack 1}-q{t+1\brack 1}-\frac{q^{n-t-1}-1}{q^{k-t-1}-1}\cdot{k-t+1\brack 1}\\
& = & \frac{q^{n-t-1}-1}{q^{k-t-1}-1}\cdot\frac{q^{k-t+1}(q^{2t-k+1}-1)}{q-1}-\frac{q(q^t-1)}{q-1}\\
& > & q^{n-k+t+1}-q^{t+1}\\
& > & 0,
\end{eqnarray*}
we obtain $h_3(t+2)>h_1(t,k+1).$

Assume that $t=\frac{k}{2}-\frac{1}{2}\leq k-3$ and $k\geq 7.$ Then $t\leq k-4.$ By (\ref{hmt+2}) and (\ref{ubound-2}), we have
$$
h_3(t+2) = {t+2\brack 1}{n-t-1\brack t}-q{t+1\brack 1}{n-t-2\brack t-1}
$$
and
$$
h_1(t,k+1) \leq {t+2\brack 1}{n-t-1\brack t}-q^{t(t-1)+1}{n-2t-2\brack t-1}{t+2\brack 2}.
$$
So by Lemma~\ref{lem1-1-1}~(iii) we have
\begin{eqnarray*}
h_3(t+2)-h_1(t,k+1) & \geq & q^{t(t-1)+1}{n-2t-2\brack t-1}{t+2\brack 2}-q{t+1\brack 1}{n-t-2\brack t-1}\\
& > & q^{(t-1)(n-2t-1)+2t+1}-q^{(t-1)(n-2t)+t+2}\\
& = & 0.
\end{eqnarray*}

Finally, assume that $t=\frac{k}{2}-\frac{1}{2}\leq k-3$ and $k=5.$ By (\ref{hmt+2}) and (\ref{hmtk+1}), we have
$$
{n-3\brack 2}={n-4\brack 1}+q^2{n-4\brack 2}
$$
and
$$
{n-2\brack 3}=\sum_{i=1}^4q^{3(i-1)}{n-2-i\brack 2}+q^{12}{n-6\brack 3}.
$$
It follows that
\begin{eqnarray*}
h_3(t+2)-h_1(t,k+1)
&=&q^3{3\brack 1}{n-4\brack 2}-\sum_{i=2}^4q^{3(i-1)}{n-2-i\brack 2}-q^4{2\brack 1}\\
&=&q^4{2\brack 1}{n-4\brack 2}-q^{6}{n-5\brack 2}-q^{9}{n-6\brack 2}-q^4{2\brack 1}\\
&=&q^4{2\brack 1}{n-5\brack 1}+q^7{n-6\brack 1}-q^4{2\brack 1}\\
&>&0.
\end{eqnarray*}

(iii)\quad This follows from the definitions of $\mathcal{H}_2(X,M,V)$ and $\mathcal{H}_3(Z)$ and the assumption that $\dim(M)=\dim(Z)=k=t+2$.   \qed

\section{Upper bounds for non-trivial $t$-intersecting families}
\label{sec:ub}

In this section we prove a number of upper bounds on the size of a maximal non-trivial $t$-intersecting family of $k$-subspaces of $V$. For any family $\mathcal{F}\subseteq{V\brack k}$ and any subspace $S$ of $V$, define
$$
\mathcal{F}_S = \{F \in \mathcal{F} \mid S \subseteq F\}.
$$

\begin{lemma}\label{lem1}
Let $\mathcal{F}\subseteq{V\brack k}$ be a $t$-intersecting family and $S$ an $s$-subspace of $V$, where $t-1\leq s\leq k-1.$ If there exists $F^\prime\in \mathcal{F}$ such that $\dim(S\cap F^\prime)=r<t,$ then for each $i\in\{1,2,\ldots,t-r\}$ there exists an $(s+i)$-subspace $T_i$ with $S\subseteq T_i$ such that $|\mathcal{F}_S|\leq {k-r\brack i}|\mathcal{F}_{T_i}|$.
\end{lemma}
\proof For each $i\in\{1,2,\ldots,t-r\}$, let
$$
\mathcal{H}_i=\left\{H\in{S+F^\prime\brack s+i}\mid S\subseteq H\right\}.
$$
Observe that $|\mathcal{H}_i|={k-r\brack i}$ by Lemma~\ref{lem4}. Since $\mathcal{F}$ is a $t$-intersecting family, for any $F\in\mathcal{F}_S$, we have $\dim(F\cap F^\prime)\geq t$, which implies that $\dim(F\cap(F^\prime+S))\geq s+t-r$ and there exists $H\in\mathcal{H}_i$ such that $H\subseteq F.$ Therefore, $\mathcal{F}_S=\cup_{H\in\mathcal{H}_i}\mathcal{F}_H.$ Let $T_i$ be a subspace in $\mathcal{H}_i$ such that $|\mathcal{F}_{H}|\leq|\mathcal{F}_{T_i}|$ for all $H\in\mathcal{H}_i.$ Then $|\mathcal{F}_S|\leq {k-r\brack i}|\mathcal{F}_{T_i}|$ as desired.  \qed

Lemma~\ref{lem4} implies that $|\mathcal{F}_T|\leq {n-\dim(T)\brack k-\dim(T)}$ for any subspace $T$ of $V$. So we have the following lemma.
\begin{lemma}\label{lem2}
Let $\mathcal{F}\subseteq{V\brack k}$ be a $t$-intersecting family and $S$ an $s$-subspace of $V$, where $t-1\leq s\leq k-1.$ If there exists $F^\prime\in \mathcal{F}$ such that $\dim(S\cap F^\prime)=r<t,$ then $|\mathcal{F}_S|\leq {k-r\brack t-r}{n-s-t+r\brack k-s-t+r}$. 
\end{lemma}

For a $t$-intersecting family $\mathcal{F}\subseteq {V\brack k}$, we define the $t$-\emph{covering number} $\tau_t(\mathcal{F})$ of $\mathcal{F}$ to be the minimum dimension of a subspace $T$ of $V$ such that $\dim(T\cap F)\geq t$ for any $F\in \mathcal{F}$. Note that for any non-trivial $t$-intersecting family $\mathcal{F}\subseteq{V\brack k}$ we have $t+1\leq \tau_t(\mathcal{F})\leq k.$

\begin{rem}
{\em	In \cite{Cao}, Cao also respectively described the structure of maximal non-trivial $t$-intersecting families with large size for finite sets and distance-regular graphs of bilinear forms by defining their $t$-covering number. It is shown that $t$-covering number is a useful notion to describe the structure of maximal non-trivial $t$-intersecting families. }
\end{rem}

\subsection{The case $\tau_t(\mathcal{F})=t+1$}

\begin{ass}
\label{hyp1}
Let $n, k$ and $t$ be positive integers with $6\leq 2k\leq n$ and $1\leq t\leq k-2$. Let $\mathcal{F}\subseteq {V\brack k}$ be a maximal non-trivial $t$-intersecting family with $\tau_t(\mathcal{F})=t+1.$ Define
$$
\mathcal{T} = \left\{T \in {V\brack t+1} \mid \dim(T\cap F)\geq t \text{ for any } F \in \mathcal{F}\right\}.
$$
\end{ass}

\begin{lemma}\label{lem3}
Let $n,\ k,\ t,\ \mathcal{F}$ and $\mathcal{T}$ be as in Assumption~\ref{hyp1}. Then $\mathcal{T}$ is a $t$-intersecting family with $t\leq\tau_t(\mathcal{T})\leq t+1$. Moreover, the following hold:
\begin{itemize}
\item[{\rm (i)}] if $\tau_t(\mathcal{T})=t$, then there exist a $t$-subspace $X$ and an $l$-subspace $M$ of $V$ with $X\subseteq M$ and $t+1\leq l\leq k+1$  such that
\begin{equation}
\label{eq:calT}
\mathcal{T} = \left\{T \in {M \brack t+1} \mid X \subseteq T\right\};
\end{equation}
\item[{\rm (ii)}] if $\tau_t(\mathcal{T})=t+1$, then there exists a $(t+2)$-subspace $Z$ of $V$ such that $\mathcal{T}={Z\brack t+1}.$
\end{itemize}
\end{lemma}
\proof The maximality of $\mathcal{F}$ implies that, for any $T\in\mathcal{T}$, $\mathcal{F}$ contains all $k$-subspaces of $V$ containing $T$. Since $2k\leq n$, for any $T_1,T_2\in \mathcal{T}$, if $\dim(T_1\cap T_2)<t,$ then there must exist $F_1,F_2\in\mathcal{F}$ such that $T_1\subseteq F_1$, $T_2\subseteq F_2$ and $\dim(F_1\cap F_2)<t$. However, this is impossible as $\mathcal{F}$ is maximal $t$-intersecting. Hence $\dim(T_1\cap T_2)\geq t$ and $\mathcal{T}\subseteq {V\brack t+1}$ is a $t$-intersecting family with $t\leq\tau_t(\mathcal{T})\leq t+1.$

(i)\quad Suppose that $\tau_t(\mathcal{T})=t$. Then there exists a $t$-subspace $X$ of $V$ such that $X$ is contained in every $(t+1)$-subspace in $\mathcal{T}.$ Assume that $M=\sum_{T\in\mathcal{T}}T$ and $\dim(M)=l.$ It suffices to prove \eqref{eq:calT} and  $t+1\leq\dim(M)\leq k+1$. Since $\tau_t(\mathcal{F})=t+1$, we have $\mathcal{F}\setminus\mathcal{F}_X\neq\emptyset.$ Let $F^\prime$ be any member of $\mathcal{F}\setminus\mathcal{F}_X.$ Observe that
 $\dim(X\cap F^\prime)\leq t-1.$ For any $T\in\mathcal{T},$ since $X\subseteq T$ and $\dim(T\cap F^\prime)\geq t,$ we have $\dim(X\cap F^\prime)=t-1$ and $\dim(T\cap (X+F^\prime))\geq t+1,$ which together imply that $\dim(X+F^\prime)=k+1$ and $T\subseteq X+F^\prime.$ Hence $M=\sum_{T\in\mathcal{T}}T\subseteq X+F^\prime$ and $t+1\leq l\leq k+1.$ It is clear that $\mathcal{T}\subseteq \left\{T\in{M\brack t+1}\mid X\subseteq T\right\}.$ Let $T^\prime$ be any $(t+1)$-subspace of $M$ with $X\subseteq T^\prime.$ For any $F\in\mathcal{F}$, if $X\subseteq F,$ then $\dim(T^\prime\cap F)\geq t$; if $X\nsubseteq F,$ then $T^\prime\subseteq X+F$ from the above discussion, which implies $\dim(T^\prime\cap F)\geq t$ by $\dim(X+F)=k+1.$ Hence $T^\prime\in\mathcal{T}$ and \eqref{eq:calT} is proved.

(ii)\quad Suppose that $\tau_t(\mathcal{T})=t+1$. Let $A,B,C\in\mathcal{T}$ be distinct subspaces such that $A\cap B$, $A\cap C$ and $B\cap C$ are pairwise distinct. Since $\mathcal{T}$ is $t$-intersecting, we have $\dim(A\cap B)=\dim(A\cap C)=\dim(B\cap C)=t,$ which together with $\dim(C)=t+1$ implies that $C=(A\cap C)+(B\cap C)\subseteq A+B$. Hence, $A+C\subseteq A+B$ and $B+C\subseteq A+B$, which imply that $A+B=A+C=B+C.$

Since $\tau_t(\mathcal{T})=t+1$, there exist three distinct subspaces $T_1,T_2,T_3\in\mathcal{T}$ such that $T_1\cap T_2$, $T_1\cap T_3$ and $T_2\cap T_3$ are pairwise distinct.  For any $T\in\mathcal{T}\setminus\{T_1,T_2,T_3\},$ if $T\cap T_1=T\cap T_2=T\cap T_3,$ then $\dim(T\cap T_1)=t,$ $T\cap T_1\subseteq T_2$ and $T\cap T_1\subseteq T_3,$ which imply that $T\cap T_1=T_1\cap T_2=T_1\cap T_3,$ a contradiction. Hence there exist $T_i, T_j\in\{T_1,T_2,T_3\}$ such that $T\cap T_i\neq T\cap T_j$ and
\[
T= (T\cap T_i)+(T\cap T_j)\subseteq T_1+T_2=T_1+T_3=T_2+T_3.
\]
Let $Z=T_1+T_2.$ Then $\mathcal{T}\subseteq{Z\brack t+1}.$ We now prove that ${Z\brack t+1}\subseteq\mathcal{T}.$ In fact, for any $F\in \mathcal{F},$ if $F\cap T_1=F\cap T_2=F\cap T_3,$ then $F\cap T_1\subseteq T_i$ for each $i\in\{1,2,3\}$. But this is impossible because $T_1\cap T_2,$ $T_1\cap T_3$ and $T_2\cap T_3$ are pairwise distinct and $\dim(F\cap T_1)\geq t$. Hence there exist $T_i,T_j\in\{T_1,T_2,T_3\}$ such that $F\cap T_i\neq F\cap T_j,$ which implies $\dim(F\cap Z)\geq t+1.$ So for any $F\in\mathcal{F}$ and $T^\prime\in{Z\brack t+1}$ we have $\dim(F\cap T^\prime)\geq t.$ Therefore, $\mathcal{T} = {Z\brack t+1}$ as desired. \qed

\begin{lemma}\label{prop3}
Let $n,\ k,\ t,\ \mathcal{F}$ and $\mathcal{T}$ be as in Assumption~\ref{hyp1}, and set  $M=\sum\limits_{T\in \mathcal{T}}T$.  Suppose that $\tau_t(\mathcal{T})=t$, $\dim(M) = k+1$ and $X$ is a $t$-subspace of $V$ which is contained in each $T\in\mathcal{T}$. Then
$$
\mathcal{F}=\left\{F\in {V\brack k}\mid X\subseteq F,\ \dim(F\cap M)\geq t+1\right\}\cup{M\brack k}.
$$
\end{lemma}
\proof It follows from the proof of Lemma~\ref{lem3} that, for any $F\in\mathcal{F}\setminus\mathcal{F}_X$, we have $M=F+X$ and hence $F\in{M\brack k}.$  Let $\mathcal{A}^\prime=\left\{F\in {V\brack k}\mid X\subseteq F,\ \dim(F\cap M)\geq t+1\right\}$ and $F^\prime$ be a fixed member of $\mathcal{F}\setminus\mathcal{F}_X$. For any $F\in\mathcal{F}_X$, we have $\dim(F\cap F^\prime)\geq t,$ $\dim(F^\prime\cap X)\leq t-1$ and $M=F^\prime+X$. Thus $\dim(F\cap M)\geq t+1$ and so $\mathcal{F}_X\subseteq \mathcal{A}^\prime.$ Note that $\mathcal{A}^\prime\cup {M\brack k}$ is a $t$-intersecting family. Therefore, $\mathcal{F}=\mathcal{A}^\prime\cup {M\brack k}$ by the maximality of $\mathcal{F}$. \qed

\begin{lemma}\label{prop3-1}
Let $n, k, t, \mathcal{F}$ and $\mathcal{T}$ be as in Assumption~\ref{hyp1}, and set $M=\sum\limits_{T\in \mathcal{T}}T$. Suppose that $\tau_t(\mathcal{T})=t$, $\dim(M) = k$ and $X$ is a $t$-subspace of $V$ which is contained in each $T\in\mathcal{T}$. Set $C=M+\sum\limits_{F\in\mathcal{F}\setminus \mathcal{F}_X}F$ and $c = \dim(C)$.  Then either $k+2\leq c\leq 2k-t$ or $c=n$. Moreover, the following hold:
\begin{itemize}
\item[{\rm(i)}] if $k+2\leq c\leq 2k-t,$ then $\mathcal{F} = \mathcal{H}_2(X, M, C)$; and
\item[{\rm(ii)}] if $c=n$, then $t\neq k-2$ and $\mathcal{F}=\mathcal{H}_2(X, M, V)$.
\end{itemize}
\end{lemma}
\proof It follows from the proof of Lemma~\ref{lem3} that, for any $F\in\mathcal{F}\setminus\mathcal{F}_X$, we have $\dim(F\cap X)=t-1$ and $M\subseteq X+F$. Since $X\subseteq M$, we then have $\dim(F\cap M)=k-1$. Note that $c \ge k+1$ by the definition of $c$.

Choose $F_1\in \mathcal{F}\setminus\mathcal{F}_X.$ Then $\dim(F_1+M)=k+1.$ If $c> k+1,$ then there exists $F_2\in \mathcal{F}\setminus\mathcal{F}_X$ such that $F_2\nsubseteq F_1+M,$ which implies $F_2\cap (F_1+M)=F_2\cap M.$ Similarly, if $c>k+2,$ then there exists $F_3\in \mathcal{F}\setminus\mathcal{F}_X$ such that $F_3\nsubseteq F_1+F_2+M,$ which implies  $F_3\cap (F_1+F_2+M)=F_3\cap M.$ Continuing, by mathematical induction we can prove  that there exist $F_1,F_2,\ldots, F_{c-k}\in\mathcal{F}\setminus\mathcal{F}_X$ such that
\begin{equation}
\label{eq:fim}
F_i \cap \left(M+\sum\limits_{j=1}^{i-1}F_j\right)=F_i\cap M
\end{equation}
for $i \in\{1,2,\ldots,c-k\}$.
If there exists $F^\prime\in\mathcal{F}$ such that $F^\prime\cap M=X$, then for any $i\in\{1,2,\ldots,c-k\}$, there exists $y_i\in F_i\setminus M$ such that $y_i\in F^\prime$ as $\dim(F^\prime \cap F_i)\geq t$ and $\dim(F^\prime \cap F_i\cap M)=t-1.$ Let $x_1,\ldots,x_t$ be a basis of $X$. By \eqref{eq:fim} and the choice of $F_1,F_2,\ldots,F_{c-k}$, one can easily show that $x_1,\ldots,x_t,y_1\ldots,y_{c-k}$ are linearly independent in $F^\prime$.

Suppose that $c\geq 2k-t+1.$ If there exists $F^\prime\in\mathcal{F}$ such that $F^\prime\cap M=X$, then by the above discussion we can obtain $c-k+t$ vectors in $F^\prime$ which are linearly independent, but this is impossible. Thus $\dim(A_1 \cap M)\geq t+1$ for any $A_1 \in \mathcal{F}_X$. By the maximality of $\mathcal{F}$, it is readily seen that any $k$-subspace $A_2$ of $V$ satisfying $\dim(A_2 \cap X)=t-1$ and $\dim(A_2 \cap M)=k-1$ must be in $\mathcal{F}$. Hence $C=V$ and $c=n.$

On the other hand, we have $c \ge k+2$, for otherwise we would have $c=k+1$ and $\dim(T\cap F)\geq t$ for any $T\in{C\brack t+1}$ with $X\subseteq T$ and any $F\in\mathcal{F}$, which imply $T\subseteq M$, a contradiction.

So far we have proved that either $k+2\leq c\leq 2k-t$ or $c=n$. It remains to prove (i) and (ii). Denote $\mathcal{A} = \mathcal{A}(X, M)$, $\mathcal{B}=\mathcal{B}(X,M,C)$ and $\mathcal{C}=\mathcal{C}(X, M, C)$.

(i)\quad Suppose that $k+2\leq c\leq 2k-t$.  Since  $\dim(F\cap X)=t-1$ and $\dim(F\cap M)=k-1$ for any $F\in\mathcal{F}\setminus\mathcal{F}_X$, we have $\mathcal{F}\setminus\mathcal{F}_X\subseteq \mathcal{C}.$  For any $F^\prime\in\mathcal{F}_X,$ if $\dim(F^\prime\cap M)\geq t+1$, then $F^\prime\in\mathcal{A};$ if $F^\prime\cap M=X,$ then $\dim(F^\prime\cap C)=c-k+t$ and so $F^\prime\in\mathcal{B}$ by the discussion above. Thus $\mathcal{F}\subseteq \mathcal{A}\cup\mathcal{B}\cup\mathcal{C}.$ It is routine to verify that $\mathcal{A}\cup\mathcal{B}\cup\mathcal{C}$ is a $t$-intersecting family. Thus, by the maximality of $\mathcal{F}$, we obtain $\mathcal{F}=\mathcal{A}\cup\mathcal{B}\cup\mathcal{C}$.

(ii)\quad Suppose that $c = n$. Then $\mathcal{F}=\mathcal{A}\cup\mathcal{C}$ by the discussion in (i) and the maximality of $\mathcal{F}$. If $t=k-2$, then $\mathcal{F}=\mathcal{H}_2(X,M,V)=\mathcal{H}_3(M)$, which implies that $\tau_t(\mathcal{F})=k,$ a contradiction. \qed

\begin{lemma}\label{prop4}
Let $n, k, t, \mathcal{F}$ and $\mathcal{T}$ be as in Assumption~\ref{hyp1}. Suppose that  $\tau_t(\mathcal{T})=t+1$ and $\mathcal{T}={Z\brack t+1}$ for some $(t+2)$-subspace $Z$ of $V$. Then $\mathcal{F}=\mathcal{H}_3(Z).$
\end{lemma}
\proof Since $\mathcal{T}={Z\brack t+1}$, we have $\dim(F\cap Z)\geq t$ for any $F\in\mathcal{F}.$ If there exists $F^\prime\in\mathcal{F}$ such that $\dim(F^\prime\cap Z)=t,$ then there exists $T^\prime\in\mathcal{T}$ such that $\dim(F^\prime\cap T^\prime)=t-1,$ a contradiction. Hence $\mathcal{F}\subseteq \mathcal{H}_3(Z).$ Since $\mathcal{H}_3(Z)$ is $t$-intersecting and $\mathcal{F}$ is maximal $t$-intersecting, we obtain $\mathcal{F}=\mathcal{H}_3(Z)$ as desired. \qed

\begin{lemma}\label{prop1}
Let $n,\ k,\ t,\ \mathcal{F}$ and $\mathcal{T}$ be as in Assumption~\ref{hyp1}.
\begin{itemize}
\item[{\rm (i)}] If $|\mathcal{T}|=1$, then
$$
|\mathcal{F}|\leq {n-t-1\brack k-t-1}+q{t+1\brack 1}{k-t\brack 1}{k-t+1\brack 1}{n-t-2\brack k-t-2}.
$$
\item[{\rm (ii)}] Suppose that $|\mathcal{T}|\geq 2$ and for some $t$-subspace $X$ and $l$-subspace $M$ of $V$ with $X\subseteq M$, $\mathcal{T}$ is a collection of $(t+1)$-subspaces of $V$ containing $X$ and contained in $M$. Then
   \begin{eqnarray}\label{equ7}
    |\mathcal{F}| & \leq & {l-t\brack 1}{n-t-1\brack k-t-1}+q^{l-t}{k-l+1\brack 1}{k-t+1\brack 1}{n-t-2\brack k-t-2}\nonumber\\
    & & +\ q^{k+1-t}{t\brack 1}{n-l\brack k-l+1}.
    \end{eqnarray}
    Moreover, if $l=t+2,$ then
   \begin{eqnarray}\label{equ8}
    |\mathcal{F}| & \leq & {l-t\brack 1}{n-t-1\brack k-t-1}+q^{l-t}{k-l+1\brack 1}{k-t+1\brack 1}{n-t-2\brack k-t-2}\nonumber\\
    & & +\ q^{2}{t\brack 1}{k-t+1\brack 1}{n-t-2\brack k-t-2}.
    \end{eqnarray}
\item[{\rm (iii)}] If $|\mathcal{T}|\geq 2$ and  $\mathcal{T}={Z\brack t+1}$ for some  $(t+2)$-subspace $Z$ of $V$, then $|\mathcal{F}|=h_3(t+2).$
\end{itemize}
\end{lemma}
\proof (i) Let $T$ be the unique $(t+1)$-subspace of $V$ in $\mathcal{T}$. Since $\dim(T\cap F)\geq t$ for any $F\in\mathcal{F}$, we have
\begin{eqnarray}\label{equ1}
\mathcal{F}=\mathcal{F}_T\ \cup \left(\bigcup_{S\in{T\brack t}}(\mathcal{F}_S\setminus\mathcal{F}_T)\right).
\end{eqnarray}

We now give an upper bound on $|\mathcal{F}_S\setminus\mathcal{F}_T|$ for any fixed $S\in{T\brack t}$. Since $\tau_t(\mathcal{F})=t+1,$ there exists $F^\prime\in\mathcal{F}\setminus \mathcal{F}_S$ such that $\dim(S\cap F^\prime)=t-1$ as $\dim(F^\prime \cap T)\geq t.$ So $T=(F^\prime \cap T)+S$ and $T\subseteq F^\prime +S.$ For any $F\in \mathcal{F}_S\setminus\mathcal{F}_T$, we have $(F\cap F^\prime)+S\subseteq F\cap (F^\prime+S).$ Since $\dim(F\cap F^\prime)\geq t$ and $\dim(F\cap F^\prime\cap S)\leq t-1$, we have $\dim(F\cap(F^\prime+S))\geq t+1.$ Hence there exists a $(t+1)$-subspace $H$ such that $H\neq T$, $S\subseteq H\subseteq S+F^\prime$ and $H\subseteq F.$ Therefore,
\begin{eqnarray}\label{equ2}
\mathcal{F}_S\setminus\mathcal{F}_T=\bigcup_{S\subseteq H\subseteq S+F^\prime, \atop H\neq T, \dim{H}=t+1} \mathcal{F}_H.
\end{eqnarray}

Consider an arbitrary $(t+1)$-subspace $H$ of $V$ satisfying $H\neq T$ and $S\subseteq H\subseteq S+F^\prime$. Since $T$ is the unique $(t+1)$-subspace of $V$ such that $\dim(T\cap F)\geq t$ for $F\in\mathcal{F}$, there exists $A\in\mathcal{F}$ such that $\dim(H\cap A)<t$. Hence $\dim(H\cap A)=t-1$ as $\dim(H\cap T)=\dim(S)=t$ and $\dim(T\cap A)\geq t.$ By Lemma~\ref{lem2}, we have $|\mathcal{F}_{H}|\leq{k-t+1\brack 1}{n-t-2\brack k-t-2}.$ By Lemma~\ref{lem4}, we obtain $|\mathcal{F}_T|\leq{n-t-1\brack k-t-1}$ and
\begin{eqnarray*}
\left|\left\{H\in{S+F^\prime\brack t+1}\mid S\subseteq H,\ H\neq T\right\}\right|={k-t+1\brack 1}-1=q{k-t\brack 1}.
\end{eqnarray*}
It follows from (\ref{equ1}) and (\ref{equ2}) that
\begin{eqnarray*}
|\mathcal{F}|\leq {n-t-1\brack k-t-1}+{t+1\brack 1}\cdot q{k-t\brack 1}{k-t+1\brack 1}{n-t-2\brack k-t-2}.
\end{eqnarray*}

(ii) We will prove the desired upper bound on $|\mathcal{F}|$ by establishing upper bounds on $|\mathcal{F}_X|$ and $|\mathcal{F}\setminus\mathcal{F}_X|$. Since $\tau_t(\mathcal{F})=t+1,$ we have $\dim(F\cap X)\geq t-1$ for any $F\in\mathcal{F}$, and there exists $F^\prime\in\mathcal{F}$ such that $\dim(X\cap F^\prime)=t-1.$ It follows from the proof of Lemma~\ref{lem3} that $X\subseteq M\subseteq X+ F^\prime.$

For any $F\in\mathcal{F}_X,$ we have $\dim(F\cap(X+ F^\prime))\geq t+1$ as $X\subseteq F$ and $\dim(F\cap F^\prime)\geq t.$ So
\begin{eqnarray}\label{equ3}
\mathcal{F}_X=\left(\bigcup_{X\subseteq H_1,\ H_1\in{M\brack t+1}}\mathcal{F}_{H_1}\right)\cup\left(\bigcup_{X\subseteq H_2,\ H_2\in{X+ F^\prime\brack t+1}\setminus{M\brack t+1}}\mathcal{F}_{H_2}\right).
\end{eqnarray}
Since by Lemma~\ref{lem4}, $|\mathcal{F}_{H_1}|\leq {n-(t+1)\brack k-(t+1)}$ for any $H_1\in{M\brack t+1}$, we have $|\bigcup_{X\subseteq H_1,\ H_1\in{M\brack t+1}}\mathcal{F}_{H_1}|\leq{l-t\brack 1}{n-(t+1)\brack k-(t+1)}.$ For any $H_2\in{X+ F^\prime\brack t+1}\setminus{M\brack t+1}$ with $X\subseteq H_2,$ we have $H_2\notin \mathcal{T}$ and so there exists $A \in\mathcal{F}$ such that $\dim(H_2\cap A)<t$. Hence $\dim(H_2\cap A)=t-1$ as $\dim(A \cap X)\geq t-1.$ It follows from Lemma~\ref{lem2} that $|\mathcal{F}_{H_2}|\leq {k-t+1\brack 1}{n-(t+1)-1\brack k-(t+1)-1}$. Note from Lemma~\ref{lem4} that
\begin{eqnarray*}
\left|\left\{H_2\in{X+ F^\prime\brack t+1}\setminus{M\brack t+1}\mid X\subseteq H_2\right\}\right|={k+1-t\brack 1}-{l-t\brack 1}=q^{l-t}{k-l+1\brack 1}.
\end{eqnarray*}
Therefore,
\begin{eqnarray}\label{equ4}
|\mathcal{F}_X|\leq {l-t\brack 1}{n-t-1\brack k-t-1}+q^{l-t}{k-l+1\brack 1}{k-t+1\brack 1}{n-t-2\brack k-t-2}.
\end{eqnarray}

For any $F\in\mathcal{F}\setminus\mathcal{F}_X$ and any $T\in\mathcal{T}$, since $\dim(F\cap X)=t-1$ and $X\nsubseteq F\cap T$, we have $T=(F\cap T)+X\subseteq F+X$. Thus, for any $F\in\mathcal{F}\setminus\mathcal{F}_X$, we have $M=\sum_{T\in\mathcal{T}}T\subseteq F+X,$ which implies $\dim(M\cap F)=l-1$. Hence $\mathcal{F}\setminus\mathcal{F}_X \subseteq \left\{F\in{V\brack k}\mid\dim(F\cap M)=l-1,\ X\nsubseteq F\right\}.$ Observe from Lemma~\ref{lem4} that the number of $k$-subspaces $F$ of $V$ satisfying $\dim(F\cap M)=l-1$ is $q^{k-l+1}{n-l\brack k-l+1}{l\brack 1}$, and the number of $k$-subspaces $F$ of $V$ satisfying $\dim(F\cap M)=l-1$ and $X\subseteq F$ is $N^\prime (t,t;k,l-1;n,n-l).$ By Lemma~\ref{lem5}, we then have
\begin{eqnarray}
|\mathcal{F}\setminus\mathcal{F}_X| & \leq & q^{k-l+1}{n-l\brack k-l+1}{l\brack 1}-N^\prime (t,t;k,l-1;n,n-l)\nonumber\\
& = & q^{k-l+1}{n-l\brack k-l+1}{l\brack 1}-q^{k-l+1}{n-l\brack k-l+1}{l-t\brack 1}\nonumber\\
& = & q^{k-t+1}{t\brack 1}{n-l\brack k-l+1}. \label{equ5}
\end{eqnarray}
Combining (\ref{equ4}) and (\ref{equ5}), we obtain \eqref{equ7}.

Now let us consider the case when $l=t+2$. From the discussion above, we have $\dim(M\cap F)=l-1=t+1$ for any $F\in\mathcal{F}\setminus\mathcal{F}_X,$ which implies  $$\mathcal{F}\setminus\mathcal{F}_X\subseteq \bigcup_{X\nsubseteq L,\ L\in{M\brack t+1}} \mathcal{F}_L.$$
For any $L\in{M\brack t+1}$ with $X\nsubseteq L$, since $L\notin\mathcal{T}$ and $\dim(F\cap M)\geq t$ for any $F\in\mathcal{F}$, there exists $F^\prime\in\mathcal{F}$ such that $\dim(F^\prime\cap L)=t-1.$ So $|\mathcal{F}_L|\leq {k-t+1\brack 1}{n-t-2\brack k-t-2}$ by Lemma~\ref{lem2}. Since by Lemma~\ref{lem4} the number of $(t+1)$-subspaces $L$ of $M$ with $X\nsubseteq L$ is equal to ${t+2\brack t+1}-{2\brack 1}$, we have
\begin{eqnarray}\label{equ6}
|\mathcal{F}\setminus\mathcal{F}_X|\leq q^{2}{t\brack 1}{k-t+1\brack 1}{n-t-2\brack k-t-2}.
\end{eqnarray}
Combining (\ref{equ4}) and (\ref{equ6}), we obtain \eqref{equ7}.

(iii) The desired equality follows from Lemma~\ref{prop4} and (\ref{hmt+2}).  \qed

\subsection{The case $\tau_t(\mathcal{F})\geq t+2$}

In \cite{AB}, Blokhuis et al. proved the following upper bound for $|\mathcal{F}|$ in the case when $t=1$.

\begin{lemma}{\rm(\cite{AB})}\label{lem11}
Let $n$ and $k$ be positive integers with $6\leq 2k\leq n$, and let $\mathcal{F}\subseteq {V\brack k}$ be a maximal intersecting family with $3\leq \tau_1(\mathcal{F})=m\leq k.$ Let $\mathcal{T}$ be the set of all $m$-subspaces $T$ of $V$ which satisfy $\dim(T\cap F)\geq 1$ for any $F\in\mathcal{F}$. Then the following hold:
\begin{itemize}
\item[{\rm (i)}] if $m=k$, then $|\mathcal{F}|\leq {k\brack 1}^k$;
\item[{\rm (ii)}] if $m<k$ and $|\mathcal{T}|\geq 2$, then
\begin{eqnarray}\label{upperb3-1-1}
|\mathcal{F}|\leq {m-1\brack 1}{k\brack 1}^{m-1}{n-m\brack k-m}+q^{2(m-1)}{k\brack 1}^{m-2}{n-m\brack k-m};
\end{eqnarray}
\item[{\rm (iii)}] if $m<k$ and $|\mathcal{T}|=1$, then
\begin{eqnarray}\label{upperb3-1-2}
|\mathcal{F}|\leq {m-1\brack 1}{m\brack 1}{k\brack 1}^{m-2}{n-m\brack k-m}+q^{m-1}{k-m+1\brack 1}{m\brack 1}{k\brack 1}^{m-1}{n-m-1\brack k-m-1}.
\end{eqnarray}
\end{itemize}
\end{lemma}

Using this lemma, we now prove the following bound for $\mathcal{F}$ with $\tau_1(\mathcal{F})< k$.

\begin{lemma}\label{prop3-1-1}
Let $n$ and $k$ be positive integers with $9\leq 2k+3\leq n$, and let $\mathcal{F}\subseteq {V\brack k}$ be a maximal intersecting family with $3\leq \tau_1(\mathcal{F}) < k.$ Then
\begin{eqnarray}\label{equ3-1}
|\mathcal{F}|\leq(q+1){k\brack 1}^2{n-3\brack k-3}+q^4{k\brack 1}{n-3\brack k-3}.
\end{eqnarray}
\end{lemma}
\proof Let $u_1(n,k,m)$ and $u_2(n,k,m)$ be the upper bounds in (\ref{upperb3-1-1}) and (\ref{upperb3-1-2}), respectively. By Lemma~\ref{lem1-1-1}(iii) and the assumption $n\geq 2k+3$, for $3\leq m<k$, we have
\begin{eqnarray*}
& & \frac{u_1(n,k,m)-u_2(n,k,m)}{{k\brack 1}^{m-2}{n-m-1\brack k-m-1}}\\
& =& \frac{q^{n-m}-1}{q^{k-m}-1}\cdot\left(q^m{m-1\brack 1}{k-m\brack 1}+q^{2(m-1)}\right)-q^{m-1}{k-m+1\brack 1}{m\brack 1}{k\brack 1}\\
& > & q^{n-k}(q^{m+k-3}+q^{2(m-1)})-q^{2k+m}\\
& > & 0.
\end{eqnarray*}
Thus, for $m\in\{3,4,\ldots,k-1\}$, we have
\begin{equation}
\label{eq:u1u2}
u_1(n,k,m)>u_2(n,k,m).
\end{equation}
For any $m\in\{3,4,\ldots,k-2\}$, by Lemma~\ref{lem1-1-1}(iii), (\ref{upperb3-1-1}) and $n\geq 2k+3$, we have
\begin{eqnarray*}
& & \frac{u_1(n,k,m)-u_1(n,k,m+1)}{{k\brack 1}^{m-2}{n-m-1\brack k-m-1}}\\
& = & \frac{q^{n-m}-1}{q^{k-m}-1}\cdot\left({m-1\brack 1}{k\brack 1}+q^{2(m-1)}\right)-{m\brack 1}{k\brack 1}^2-q^{2m}{k\brack 1}\\
& > & q^{n-k}(q^{m+k-3}+q^{2(m-1)})-q^{2k+m}-q^{2m+k}\\
& > & 0.
\end{eqnarray*}
So $u_1(n,k,m)$ is decreasing as $m \in \{3,4,\ldots,k-1\}$ increases. Combining this with \eqref{eq:u1u2} and Lemma \ref{lem11}, we obtain $|\mathcal{F}| \leq u_1(n, 3, m)$ for $m = \tau_1(\mathcal{F})$, which yields (\ref{equ3-1}) as $u_1(n, 3, m)$ is exactly the right-hand side of \eqref{equ3-1}. \qed

\begin{lemma}\label{prop2}
Let $n, k$ and $t$ be positive integers with $8\leq 2k\leq n$ and $2\leq t\leq k-2$, and let $\mathcal{F}\subseteq {V\brack k}$ be a maximal $t$-intersecting family with $t+2\leq \tau_t(\mathcal{F})=m\leq k.$  Then 
\begin{equation*}
|\mathcal{F}|\leq {m\brack t} {k\brack 1}^{m-t-2}{k-t+1\brack 1}^2{n-m\brack k-m}.
\end{equation*}
Moreover, if $n\geq 2k+t+1$, then
\begin{equation*}
|\mathcal{F}|\leq {t+2\brack 2} {k-t+1\brack 1}^2{n-t-2\brack k-t-2}.
\end{equation*}
\end{lemma}
\proof
Let $T$ be an $m$-subspace of $V$ which satisfies $\dim(T\cap F)\geq t$ for any $F\in\mathcal{F}$. Then $\mathcal{F}=\cup_{H\in{T\brack t}}\mathcal{F}_H$ and hence there exists $H_1\in{T\brack t}$ such that $|\mathcal{F}|\leq {m\brack t}|\mathcal{F}_{H_1}|.$ 
If $m\geq t+3,$ using Lemma~\ref{lem1} repeatedly, then there exist $H_2\in{V\brack t+1}$, $H_3\in{V\brack t+2}$,\ldots, $H_{m-t-1}\in{V\brack m-2}$ such that $H_i\subseteq H_{i+1}$ and $|\mathcal{F}_{H_i}|\leq{k\brack 1}|\mathcal{F}_{H_{i+1}}|$ for each $i\in\{1,2,\ldots,m-t-2\}$. Thus there exists $H^\prime\in{V\brack m-2}$ such that
\begin{equation*}
|\mathcal{F}|\leq{m\brack t}{k\brack 1}^{m-t-2}|\mathcal{F}_{H^\prime}|.
\end{equation*}
Since $\tau_t(\mathcal{F})> m-2,$ we have $\mathcal{F}\setminus\mathcal{F}_{H^\prime}\neq\emptyset$ and $\dim(F\cap H^\prime)\leq t-1$ for any $F\in\mathcal{F}\setminus\mathcal{F}_{H^\prime}.$

\medskip
\textsf{Case 1.}\quad $\dim(F\cap H^\prime)\leq t-2$ for all $F\in\mathcal{F}\setminus\mathcal{F}_{H^\prime}.$
\medskip

Let $F_1$ be a fixed  $k$-subspace in $\mathcal{F}\setminus\mathcal{F}_{H^\prime}.$  Let $s_1 = \dim(F_1\cap H^\prime)$ so that $0\leq s_1\leq t-2$. By Lemma~\ref{lem2}, we have
\begin{equation*}
|\mathcal{F}_{H^\prime}|\leq{k-s_1\brack t-s_1}{n-m+2-t+s_1\brack k-m+2-t+s_1},
\end{equation*}
which implies that
\begin{equation}
\label{eq:calF}
|\mathcal{F}|\leq{m\brack t}{k\brack 1}^{m-t-2}{k-s_1\brack t-s_1}{n-m+2-t+s_1\brack k-m+2-t+s_1}.
\end{equation}
Let
$$
g(s)={k-s\brack t-s}{n-m+2-t+s\brack k-m+2-t+s}
$$
for $s \in\{0,1,\ldots,t-2\}$. By $n\geq 2k$ and Lemma~\ref{lem1-1-1}(ii), we have
$$
\frac{g(s+1)}{g(s)}=\frac{(q^{t-s}-1)(q^{n-m+3-t+s}-1)}{(q^{k-s}-1)(q^{k-m+3-t+s}-1)}>q^{n-2k+t-1}> 1
$$
for $s\in\{0,1,\ldots,t-3\}$. That is, the function $g(s)$ is increasing as $s \in\{0,1,\ldots,t-2\}$ increases. This together with \eqref{eq:calF} yields
\begin{equation}
\label{equ10}
|\mathcal{F}| \le {m\brack t}g(s_1) \le {m\brack t}g(t-2) = {m\brack t}{k\brack 1}^{m-t-2}{k-t+2\brack 2}{n-t-2\brack k-t-2}.
\end{equation}

\medskip
\textsf{Case 2.}\quad There exists $F_2\in\mathcal{F}\setminus\mathcal{F}_{H^\prime}$ such that $\dim(F_2\cap H^\prime)=t-1.$
\medskip

By Lemma~\ref{lem1}, there exists an $(m-1)$-subspace $H^{\prime\prime}$ such that $|\mathcal{F}_{H^\prime}|\leq{k-t+1\brack 1}|\mathcal{F}_{H^{\prime\prime}}|$. Hence $|\mathcal{F}|\leq {m\brack t}{k\brack 1}^{m-t-2}{k-t+1\brack 1}|\mathcal{F}_{H^{\prime\prime}}|.$
Since $\tau_t(\mathcal{F})>m-1$, there exists $F_3\in \mathcal{F}$ such that $\dim(F_3\cap H^{\prime\prime})\leq t-1.$ 

If $\dim(F_3\cap H^{\prime\prime})=t-1$, then there exists an $m$-subspace $H^{\prime\prime\prime}$ with $H^{\prime\prime}\subseteq H^{\prime\prime\prime}$ such that $|\mathcal{F}_{H^{\prime\prime}}|\leq{k-t+1\brack 1}|\mathcal{F}_{H^{\prime\prime\prime}}|.$ Since $|\mathcal{F}_{H^{\prime\prime\prime}}|\leq{n-m\brack k-m}$ by Lemma~\ref{lem4}, we have
\begin{eqnarray}\label{equ11}
|\mathcal{F}|\leq{m\brack t}{k\brack 1}^{m-t-2}{k-t+1\brack 1}^2{n-m\brack k-m}.
\end{eqnarray}

Suppose that $\dim(F_3\cap H^{\prime\prime})=s_2\leq t-2$. by Lemma~\ref{lem2}, we have
$$
|\mathcal{F}_{H^{\prime\prime}}|\leq{k-s_2\brack t-s_2}{n-m+1-t+s_2\brack k-m+1-t+s_2}.
$$
Similar to Case 1, it is straightforward to verify that the function ${k-s\brack t-s}{n-m+1-t+s\brack k-m+1-t+s}$ is increasing as $s\in\{0,1,\ldots,t-2\}$ increases. Hence
\begin{eqnarray}\label{equ11-1}
|\mathcal{F}|\leq{m\brack t}{k\brack 1}^{m-t-2}{k-t+1\brack 1}{k-t+2\brack 2}{n-m-1\brack k-m-1}.
\end{eqnarray}
\medskip

By Lemma~\ref{lem1-1-1}~(ii) and $n\geq 2k$, it is straightforward to verify that
$$
{k-t+1\brack 1}^2{n-m\brack k-m}\geq\max\left\{{k-t+2\brack 2}{n-m\brack k-m},\ {k-t+1\brack 1}{k-t+2\brack 2}{n-m-1\brack k-m-1}\right\}.
$$
This together with \eqref{equ10}, \eqref{equ11} and \eqref{equ11-1} yields
\begin{equation}
\label{eq:eq1}
|\mathcal{F}|\leq{m\brack t}{k\brack 1}^{m-t-2}{k-t+1\brack 1}^2{n-m\brack k-m}.
\end{equation}

Let 
$$
p(m^\prime)={m^\prime\brack t}{k\brack 1}^{m^\prime-t-2}{n-m^\prime\brack k-m^\prime}
$$
for $m^\prime\in\{t+2,t+3,\ldots,k\}$. By $n\geq 2k$ and Lemma~\ref{lem1-1-1}~(ii), we have
$$
\frac{p(m^\prime)}{p(m^\prime+1)}=\frac{(q^{m^\prime-t+1}-1)(q-1)(q^{n-m^\prime}-1)}{(q^{m^\prime+1}-1)(q^k-1)(q^{k-m^\prime}-1)}>q^{n-2k-t-1}\geq 1
$$
for $m^\prime\in\{t+2,t+3,\ldots,k-1\}$. That is, the function $q^(m^\prime)$ is decreasing as $m^\prime\in \{t+2,t+3,\ldots,k\}$ increases. This together with \eqref{eq:eq1} yields
$$
|\mathcal{F}|\leq {k-t+1\brack 1}^2p(t+2)\leq{t+2\brack 2} {k-t+1\brack 1}^2{n-t-2\brack k-t-2}.
$$
Therefore, the desired upper bounds follow.  \qed

\section{Proofs of Theorems \ref{main-1} and \ref{main-2}}
\label{sec:pf}

\noindent\textit{Proof of Theorem~\ref{main-1}.}\quad
Let $n, k$ and $t$ be positive integers with $k\geq 3$. Suppose that $n \ge 2k+t+\min\{4, 2t\}$. That is, if $t=1$, then $n \ge 2k+3$, and if $t \ge 2$, then $n \ge 2k+t+4$. Let $\mathcal{F}\subseteq {V\brack k}$ be a maximal non-trivial $t$-intersecting family which is not any of the exceptional families in (i) and (ii) of Theorem~\ref{main-1}. Set
$$
f_2(n,k,t) = \frac{f(n,k,t)-|\mathcal{F}|}{{n-t-2\brack k-t-2}},
$$
where the function $f$ is as defined in \eqref{eq:f}. It suffices to prove  $f(n,k,t)>|\mathcal{F}|$ or equivalently $f_2(n,k,t)>0$.

Let $\mathcal{T}$ be the set of all $\tau_t(\mathcal{F})$-subspaces $T$ of $V$ which satisfy $\dim(T\cap F)\geq t$ for any $F\in\mathcal{F}$.

\medskip
\textsf{Case 1.}\quad $\tau_t(\mathcal{F})=t+1.$

\medskip
\textsf{Case 1.1.}\quad $|\mathcal{T}|=1$.
\medskip

In this case, by Lemma~\ref{prop1}(i), we have
$$
q^{-1} f_2(n,k,t)\geq{n-t-1\brack 1}-{k-t\brack 2}-{t+1\brack 1}{k-t\brack 1}{k-t+1\brack 1}.
$$
If $(t,q)=(1,2)$, then $n\geq 2k+3$ and
\begin{eqnarray*}
\frac{3}{2}\cdot f_2(n,k,t) & \geq & 3\cdot (2^{n-2}-1)-(2^{k-1}-1)(2^{k-2}-1)-9\cdot (2^{k-1}-1)(2^{k}-1)\\
& = & 3\cdot 2^{n-2}-37\cdot 2^{2k-3}+9\cdot 2^k+10\cdot 2^{k-1}+2^{k-2}-13\\
& > & 0
\end{eqnarray*}
as desired.
Suppose that $n \ge 2k+3$ and $q\geq3$ if $t=1$, and that $n \geq 2k+t+4$ if $t\geq 2.$ By Lemma~\ref{lem1-1-1}(iii)(iv), we have
$$
q^{-1} f_2(n,k,t)>q^{n-t-2}-q^{2(k-t-1)}-8q^{2k-t-1}>0
$$
as desired.

\medskip
\textsf{Case 1.2.}\quad $|\mathcal{T}|\geq 2$ and $\tau_t(\mathcal{T})=t.$
\medskip

Let $M=\sum_{T\in \mathcal{T}} T$ and $l=\dim(M)$. Since $\mathcal{F}$ is a maximal non-trivial $t$-intersecting family other than any of the exceptional families in Theorem~\ref{main-1}, we have $l\leq k-1$ by Lemmas~\ref{prop3} and \ref{prop3-1}.

Let us first consider the case when $l=t+2$. Then $k\geq 4$ as $l\leq k-1$. By (\ref{equ8}) and Lemma~\ref{lem1-1-1}(iii), we have
\begin{eqnarray*}
& & q^{-2} f_2(n,k,t)\\
& \geq & \frac{q^{n-t-1}-1}{q^{k-t-1}-1} {k-t-2\brack 1}-q^{-1}{k-t\brack 2}-{k-t-1\brack 1}{k-t+1\brack 1} -{t\brack 1}{k-t+1\brack 1}\\
& > & q^{n-t-3}-q^{2k-2t-3}-q^{2k-2t}-q^{k+1}.
\end{eqnarray*}
Thus, if $t=1$, then
$$
q^{-2} f_2(n,k,t)>q^{k+1}\left(q^{n-k-5}-q^{k-6}-q^{k-3}-1\right)>0
$$
as $n\geq 2k+3$. If $t\geq 2$, then
$$
n-t-3\geq \max\{2k-2t,\ k+1\}+2
$$
as $n \geq 2k+t+4$, and hence the inequality above implies $q^{-2} f_2(n,k,t)>0.$ In either case we have $f_2(n,k,t)>0$ as desired.

Now consider the case when $t+3\leq l\leq k-1.$ Then
$$
{n-l\brack k-l+1}\leq {n-t-3\brack k-t-2}.
$$
Since $n-t-1\geq \max\{2k-2t+2,\ t+3\}+2$ and $t\leq k-2$, by (\ref{equ7}) we have
\begin{eqnarray*}
& & f_2(n,k,t)\\
& \geq & \frac{q^{n-t-1}-1}{q^{k-t-1}-1} q^{l-t} \cdot {k-l\brack 1}-q{k-t\brack 2}-q^{l-t}{k-l+1\brack 1}{k-t+1\brack 1}
    -q^{k+1-t}{t\brack 1}\cdot \frac{{n-l\brack k-l+1}}{{n-t-2\brack k-t-2}}\\
& \geq & \frac{q^{n-t-1}-1}{q^{k-t-1}-1} q^{l-t}\cdot{k-l\brack 1}-q{k-t\brack 2}-q^{l-t}{k-l+1\brack 1}{k-t+1\brack 1}
    -q^{k+1-t}{t\brack 1} \cdot \frac{q^{n-k}-1}{q^{n-t-2}-1}\\
& > &q^{n-t-1}-q^{2k-2t-1}-q^{2k-2t+2}-q^{t+3}\\
& > & 0,
\end{eqnarray*}
as desired.

\medskip
\textsf{Case 1.3.}\quad $|\mathcal{T}|\geq 2$ and $\tau_t(\mathcal{T})=t+1.$
\medskip

In this case, by Lemma~\ref{prop1}(i) and Lemma~\ref{lem10}(i), we have
$f(n,k,t)>|\mathcal{F}|$ if $1\leq t\leq \frac{k}{2}-\frac{3}{2}.$

\medskip
\textsf{Case 2.}\quad $t+2\leq\tau_t(\mathcal{F})\leq k.$

\medskip
\textsf{Case 2.1.}\quad $t=1.$
\medskip

Since $t=1$, we have $3 \le \tau_1(\mathcal{F}) \le k$. Consider the case $\tau_1(\mathcal{F})=k$ first. By Lemma~\ref{lem11}(i), in this case we have
\begin{equation}\label{equ17}
f(n,k,t)-|\mathcal{F}| \geq {k-1\brack 1}{n-2\brack k-2}-q{k-1\brack 2}{n-3\brack k-3}-{k\brack 1}^k.
\end{equation}
If $(n,k,q)=(9,3,2)$, then $f(n,k,1)-|\mathcal{F}|\geq 36>0.$ If $(n,k)=(9,3)$ and $q\geq3$, then $f(n,k,1)-|\mathcal{F}|\geq q^7+q^6-q^5-4q^4-5q^3-4q^2-2q>0$. Since when $k=3$ the right-hand side of (\ref{equ17}) is increasing with $n$, we have $f(n,3,1)-|\mathcal{F}|>0$ for $n\geq 10.$ If $(n,k,q)=(11,4,2)$, then $f(n,k,1)-|\mathcal{F}|\geq 249850>0.$ If $(n,k)=(11,4)$ and $q\geq 3$, or $n=2k+3$ and $k\geq5$, or $n\geq 2k+4$ and $k\geq 4$, then by Lemma~\ref{lem1-1-1}(iii)(iv),
\begin{eqnarray*}
f(n,k,1)-|\mathcal{F}| & > & {n-3\brack k-3}(q^{n-2}-q^{2k-3})-2^k\cdot q^{k(k-1)}\\
& > & q^{(k-3)(n-k)+n-3}-2^k\cdot q^{k(k-1)}\\
& \geq & 0.
\end{eqnarray*}

Now let us consider the case when $3\leq \tau_1(\mathcal{F})<k$. In this case, by Lemma~\ref{prop3-1-1}, we have
$$
  f_2(n,k,1)\geq\frac{q^{n-2}-1}{q^{k-2}-1} {k-1\brack 1}-q{k-1\brack 2}-(q+1){k\brack 1}^2 -q^4{k\brack 1}.
$$
If $q\geq 3,$ then by $n\geq 2k+3$, $k\geq 4$ and Lemma~\ref{lem1-1-1}(iii)(iv), we obtain
\begin{eqnarray*}
f_2(n,k,1) & > & q^{n-2}-q^{2k-3}-4(q+1)q^{2k-2}-q^{k+4}\\
& \geq & q^{2k-3}\left(q^4-1-4q^2-4q-q^{7-k}\right)\\
& > & 0.
\end{eqnarray*}
If $q=2$, then by $n\geq 2k+3$, $k\geq 4$ and $(2^{n-2}-1)(2^{k-1}-1)/(2^{k-2}-1)>2^{2k+2}$ we obtain
\begin{eqnarray*}
  f_2(n,k,1)& \geq & \frac{(2^{n-2}-1)(2^{k-1}-1)}{2^{k-2}-1}-\frac{2}{3}(2^{k-1}-1)(2^{k-2}-1)-3(2^k-1)^2
-2^4 (2^k-1)\\
& > & \frac{11}{3} \cdot 2^{2k-2}-19\cdot 2^{k-1}+\frac{37}{3}\\
& > & 0.
\end{eqnarray*}

\textsf{Case 2.2.}\quad $t\geq 2.$
\medskip

By Lemma~\ref{prop2}, we have 
$$
f_2(n,k,t)\geq\frac{q^{n-t-1}-1}{q^{k-t-1}-1}{k-t\brack 1}-q{k-t\brack 2}-{t+2\brack 2}{k-t+1\brack 1}^2.
$$

Assume that $n=2k+t+4$ and $q=2$ first. We have
\begin{align*}
&f_2(n,k,t)\\
=&\frac{(2^{2k+3}-1)(2^{k-t}-1)}{2^{k-t-1}-1}-\frac{2}{3}(2^{k-t}-1)(2^{k-t-1}-1)-\frac{1}{3}(2^{t+2}-1)(2^{t+1}-1)(2^{k-t+1}-1)^2\\
>&2^{2k+4}-\frac{1}{3}\cdot 2^{2k-2t}-\frac{1}{3}\cdot 2^{2k+5}\\
>&0.
\end{align*}

Now assume that $n=2k+t+4$ and $q\geq3$, or $n\geq 2k+t+5$. Since $q{k-t\brack 2}<{k-t+1\brack1}^2$ and ${t+2\brack 2}+1\leq4q^{2t},$ by Lemma~\ref{lem1-1-1}~(ii)(iv), we have
\begin{align*}
f_2(n,k,t)&>\frac{q^{n-t-1}-1}{q^{k-t-1}-1}{k-t\brack 1}-4q^{2t}{k-t+1\brack 1}^2\\
&={k-t\brack 1}\left(\frac{q^{n-t-1}-1}{q^{k-t-1}-1}-\frac{4q^{2t}(q^{k-t+1}-1)^2}{(q-1)(q^{k-t}-1)}\right)\\
&>{k-t\brack 1}\left(q^{n-k}-16q^{k+t+1}\right)\\
&\geq 0.
\end{align*}
This completes the proof.
\qed

\medskip
\noindent\textit{Proof of Theorem~\ref{main-2}.}\quad The result follows from Theorem~\ref{main-1}, Remark~\ref{rem1} and Lemmas~\ref{lem6}, \ref{lem7}, \ref{lem10} and \ref{lem8-1}. \qed

\section*{Acknowledgement}

This research was supported by NSFC (11671043) and NSF of Hebei Province (A2019205092).

%\end{CJK*}

\end{document}